\documentclass{amsart}
\usepackage{amsmath,amsfonts,amsthm,amssymb,mathrsfs,latexsym,paralist}
\usepackage{enumerate}
\usepackage{nomencl}
\usepackage{tikz, soul}
\usepackage{dsfont}
\usepackage{mathabx}

\newtheorem{theorem}{Theorem}[section]
\newtheorem{lemma}[theorem]{Lemma}
\newtheorem{proposition}[theorem]{Proposition}

\theoremstyle{definition}
\newtheorem{definition}[theorem]{Definition}
\newtheorem{remark}[theorem]{Remark}

\numberwithin{equation}{section}

\newcommand{\co}{\mskip0.5mu\colon\thinspace}
\newcommand{\ve}{\varepsilon}
\newcommand{\C}{\mathbb{C}}

\newcommand{\N}{\mathbb{N}}

\newcommand{\R}{\mathbb{R}}

\newcommand{\Z}{\mathbb{Z}}

\newcommand{\HH}{\mathcal{H}}

\newcommand{\HHDkK}{\mathcal{H}_D^p(K)}

\newcommand{\LOkK}{L^2\Omega^p(K)}

\newcommand{\HOKE}{{H}^1\Omega_D^p(K)}

\newcommand{\HHDkKo}{\mathcal{H}_D^p(K)^{\obot}}
\newcommand{\HHNkKo}{\mathcal{H}_N^p(K)^{\obot}}

\newcommand{\DP}{\varphi_{D}}
\newcommand{\bft}{\mathbf{t}}
\newcommand{\bfn}{\mathbf{n}}

\newcommand{\PK}{\partial {K}}
\newcommand{\FLOW}{\psi^X_t}
\newcommand{\PB}{(\psi^X_t)^*}
\newcommand{\SIGMA}{\Sigma_t^X}

\begin{document}

\title{
Approximating orbifold spectra using collapsing connected sums
}

\author{Carla Farsi}
\address{Department of Mathematics, University of Colorado at Boulder,
    UCB 395, Boulder, CO 80309-0395}
\email{farsi@euclid.colorado.edu}

\author{Emily Proctor}
\address{Department of Mathematics, Middlebury College, Middlebury, VT 05753}
\email{eproctor@middlebury.edu}

\author{Christopher Seaton}
\address{Department of Mathematics and Computer Science,
    Rhodes College, 2000 N. Parkway, Memphis, TN 38112}
\email{seatonc@rhodes.edu}



\begin{abstract}
For a closed Riemannian orbifold $O$, we compare the spectra of the Laplacian, acting on functions or differential forms, to the Neumann spectra of the orbifold with boundary given by a domain $U$ in $O$ whose boundary is a smooth manifold. Generalizing results of several authors, we prove that the metric of $O$ can be perturbed to ensure that the first $N$ eigenvalues of $U$ and $O$ are arbitrarily close to one another. This involves a generalization of the Hodge decomposition to the case of orbifolds with manifold boundary. Using these results, we study the behavior of the Laplace spectrum on functions or forms of a connected sum of two Riemannian orbifolds as one orbifold in the pair is collapsed to a point. We show that the limits of the eigenvalues of the connected sum are equal to those of the non-collapsed orbifold in the pair. In doing so, we prove the existence of a sequence of orbifolds with singular points whose eigenvalue spectra come arbitrarily close to the spectrum of a manifold, and a sequence of manifolds whose eigenvalue spectra come arbitrarily close to the eigenvalue spectrum of an orbifold with singular points. We also consider the question of prescribing the first part of the spectrum of an orientable orbifold.
\end{abstract}

\maketitle

\tableofcontents


\section{Introduction}
\label{sec:intro}

In a series of articles, Ann\'e, Colbois, Takahashi, and others have studied the behavior of eigenvalues of the Laplacian acting on functions and forms on closed orientable manifolds with no curvature assumptions under conditions of collapsing a submanifold \cite{anneperturbNeumann,annespectre,annecolboisoperateur,annecolbois,Anne-Post,annetakahashi,colboiselsufi,takahashi}.
These authors have considered handled manifolds, manifolds with balls removed, and connected sums of manifolds.  They showed that as the handles collapse, the balls fill in, or one part of the connected sum collapses, the limit of the spectrum of the Laplacian acting on functions (resp. $p$-forms) is equal to the $0$-spectrum (resp. $p$-spectrum) of the Laplace spectrum of the limit space, with careful counting of the zeros of the spectrum. Note that similar results about the convergence of eigenvalues have been established for manifolds with Ricci curvature bounded below, see \cite{Wei} for a survey.

In this paper, we consider generalizations of some of these results to the case of orbifolds. We focus on the case of the connected sum $O_1 \sqcup O_2$ of two orbifolds $O_1$ and $O_2$ along a nonsingular locus and the behavior of the spectrum of the Laplacian acting on functions or $p$-forms as one of these orbifolds is collapsed to a point. When $O_1$ and $O_2$ are manifolds, Takahashi proved in \cite{takahashi} that when $O_2$ collapses to a point, the $0$-spectrum of $O_1\sqcup O_2$
converges to that of $O_1$. An analogous result for the spectrum of the Laplacian acting on $p$-forms was later proven by Ann\'{e} and Takahashi in \cite{annetakahashi}.

Our primary motivation for considering this generalization is related to the well-studied yet open question of whether the spectrum of the Laplacian acting on functions or forms can detect the existence of orbifold singularities. In other words, it is unknown whether or not there can be a manifold that is isospectral to an orbifold with singular points. Doyle and Rossetti (\cite{doylerossetti}) and Rossetti, Schueth, and Weilandt (\cite{rsw}) have produced examples of isospectral pairs of orbifolds having different maximal isotropy orders.  This indicates at least the possibility of an isospectral manifold-orbifold pair.  We note that in the same direction, Gordon and Rossetti showed in \cite{gordonrossetti} that it is possible to have a $2p$-dimensional manifold that is isospectral on $p$-forms to an orbifold with a nontrivial singular set.

In certain contexts, however, it has been shown that a manifold cannot be $0$-isospectral to an orbifold with singular points.  In \cite{drydenstrohmaier}, Dryden and Strohmaier showed that isospectral orientable hyperbolic orbisurfaces have the same number of cone points of a particular order.   Linowitz and Meyer have also shown that under mild assumptions, there is no such manifold-orbifold pair in the class of length-commensurable compact locally symmetric spaces of nonpositive curvature associated to simple Lie groups (\cite{linowitzmeyer}).  More generally, by the work of Dryden, Gordon, Greenwald and Webb it is impossible for an even- (resp. odd-) dimensional orbifold having an odd- (resp. even-) dimensional singular stratum to be isospectral to a manifold (\cite{dggw}).  In \cite{gordonrossetti}, Gordon and Rossetti showed that a manifold and an orbifold with singular points having a common Riemannian covering cannot be isospectral. Finally, Sutton has shown that it is even the case that if an orbifold $O$ with singular points and a manifold $M$ admit isospectral Riemannian coverings $M_1$ and $M_2$, then $O$ and $M$ cannot be isospectral (\cite{sutton}).

Here, we extend Takahashi and Ann\'{e}--Takahashi's results to demonstrate that if the orbifold $O_2$ in a connected sum
$O_1\sqcup O_2$ is collapsed to a point, then the spectrum of the Laplacian acting on functions or forms converges to the corresponding spectrum of $O_1$. To accomplish this, we generalize to orbifolds the results of Colin de Verdi\`ere in
\cite{VerdierMultiplicite}, demonstrating that on a closed orbifold $O$ of dimension at least $3$, the metric can be
perturbed off a domain $U$ so that the first $N$ eigenvalues of the Laplacian acting on functions on $O$ are arbitrarily close to the eigenvalues of the Laplacian acting on functions on $\overline{U}$ with Neumann boundary conditions, see
Theorem~\ref{thrm:SpecDefDim3}. This formulation is framed in terms of the so-called
$N$-spectral defect of the quadratic forms associated to the respective operators, roughly the maximum difference
between the first $N$ corresponding eigenvalues, see Definition~\ref{def:N-spectral-defect}. Because our intended application is that of connected sums, we restrict to the case of orbifold domains $U$ with smooth boundary and no singular points on the boundary. We also describe the extension to orbifolds of a similar result of Rauch and Taylor \cite{rauchtaylor}, extended to more general contexts in \cite{anneperturbNeumann,annespectre,annecolboisoperateur,Anne-Post,chavel-feldman-domains,chavel-feldman-less},
comparing the eigenvalues on a Euclidean domain to those of the Neumann boundary conditions after removing a small
ball, see Theorem~\ref{thrm:RauchTaylor}.
For the analogous result for the Laplacian acting on differential forms, we generalize a result of Jammes
\cite{JammesPrecrip} to orbifolds with manifold boundary, see Theorem~\ref{thrm:SpecDefForms}. An essential
ingredient, which may be of independent interest, is an extension of the Hodge decomposition to such orbifolds; we explain this extension in Theorem~\ref{thrm:Hodge} following the exposition of \cite{schwarz}. Note that the Hodge decomposition has been extended to closed orbifolds in \cite{bailey}, but an extension to the case of orbifolds with boundary is missing from
the literature.

We then apply these results to collapsing $O_2$ in a connected sum $O_1\sqcup O_2$ of orbifolds
admitting a piecewise differentiable metric that is not differentiable along the boundary identification as in Takahashi \cite{takahashi}.
For the Laplacian acting on functions, we give two formulations in
Theorems~\ref{thrm:RTsmoothconvergence} and \ref{thrm:smoothconvergence} that differ in whether one allows
perturbation of the original metric on $O_1$; the result for the Laplacian acting on differential forms is
given in Theorem~\ref{thrm:smoothconvergencePForms}.
By choosing one of the $O_i$ to be a smooth manifold and approximating the piecewise differentiable metric on the connected sum with a smooth orbifold metric, this yields a sequence of orbifolds with singular points whose spectra (on functions or forms) converge to that of a smooth manifold, and similarly a sequence of smooth manifolds whose spectra converge to that of an orbifold with singular points; see Theorems~\ref{thrm:orbifoldstomanifold}, \ref{thrm:manifoldstoorbifold}, \ref{thrm:orbifoldstomanifoldForms}, and \ref{thrm:manifoldstoorbifoldForms}. Hence, while the results of this paper do not provide the existence of a manifold that is isospectral to an orbifold with singular points, we demonstrate that the spectra of a nontrivial orbifold and smooth manifold can in some sense be arbitrarily near one another using the intuitive construction of collapsing a connected sum.
We also extend results of Colin de Verdi\`ere by showing that the first $N$ eigenvalues of the Laplacian acting on functions on a closed orbifold $O$ of dimension $n\geq 3$ can be prescribed.
In particular, given a closed Riemannian manifold $(M, g^\prime)$ of any dimension and any closed oriented orbifold of dimension $\geq 3$, Theorem~\ref{thrm:PrescribNSpec} demonstrates that one can find a metric $g$ on $O$ such that the first $N$ nonzero eigenvalues of the Laplacian of $(O,g)$ acting on functions are equal to the first $N$ nonzero eigenvalues of the Laplacian of $(M,g^\prime)$ on functions. Switching the roles of $M$ and $O$, already by the results of \cite{VerdierConstrucDonee}, a manifold admits a metric with respect to which its first $N$ eigenvalues coincide with those of any fixed Riemannian orbifold (and similarly for the Laplacian acting on forms using the results of \cite{JammesPrecrip}), though the relationship between $M$ and $O$ does not have the clear intuitive connection of collapsing of connected sums.

The outline of this paper is as follows. Sections~\ref{sec:OfldPrelim} through \ref{sec:Laplacian} review background and well-established results on manifolds, describing the modifications required to extend them to the orbifold case under consideration. The reader familiar with orbifolds may begin reading Sections~\ref{sec:SpecDefDim3Func} and refer back to these earlier sections for definitions and background. Sections~\ref{sec:SpecDefDim3Func} through \ref{sec:ConnectedSums} include the key results on collapsing. Specifically, in Section~\ref{sec:OfldPrelim}, we give background on orbifolds with boundary and define the piecewise differentiable metrics that will appear on connected sum orbifolds. Section~\ref{sec:SobolevSpaces} explains the extensions of fundamental results on differential forms such as Rellich's theorem, the Trace theorem, and Green's theorem to orbifolds with boundary as well as to the case of piecewise differentiable connected sum metrics. The extension of the Hodge decomposition theorem to orbifolds with manifold boundary is given in Section~\ref{sec:Hodge}, and Section~\ref{sec:Laplacian} confirms familiar properties of the Laplacian acting on $p$-forms in the case of orbifolds with boundary (using various boundary conditions) and orbifold connected sums. In Sections~\ref{sec:SpecDefDim3Func} and \ref{sec:SpecDefDim3Form}, for the Laplacian acting on functions and differential forms, we compare the spectrum of a closed orbifold to the spectrum of an orbifold domain with manifold boundary and Neumann boundary conditions. Section~\ref{sec:ConnectedSums} gives the applications of these results to collapsing connected sums and the question of prescribing the first $N$ eigenvalues of the Laplacian acting on functions for a closed orbifold $O$.


\section*{Acknowledgements}

The authors thank Colette Ann\'e and Junya Takahashi for helpful correspondence in the course of this work.
The authors would also like to thank the anonymous referees for their very useful comments and suggestions that significantly improved this paper.

C.S. and E.P. would like to thank the Department of Mathematics at the University of Colorado at Boulder, and C.F. and E.P. would like to thank the Department of Mathematics and Computer Science at Rhodes College, for hospitality during work on this manuscript.
C.F. would like to thank the sabbatical program at the University of Colorado at Boulder and was partially supported by the Simons Foundation Collaboration Grant for Mathematicians \#523991.
C.S. would like to thank the sabbatical program at Rhodes College and was partially supported by the E.C. Ellett Professorship in Mathematics.


\section{Orbifolds with boundary and with connected sum metrics}
\label{sec:OfldPrelim}

In this section, we introduce definitions of our two main objects of study, orbifolds with boundary and orbifolds with connected sum metrics.


\subsection{Orbifolds with boundary}
\label{subsec:OfldsBoundary}

In order to carefully define an orbifold with boundary, we begin by defining a real half space.  Our treatment here follows that of \cite{alr}, \cite{chiang}, and \cite{schwarz}.

Given a nonzero vector $u\in\mathbb{R}^n$, the corresponding half space $\R^n_u$ is given by
\[
\R^n_u:=\{x\in\R^n \,\vert\, \langle x,u\rangle\geq 0\}.
\]
We say that a map $h\co \R^n_u\to \R^k$ is differentiable at a boundary point of $\R^n_u$ if it has a differentiable extension $\tilde{h}\co\R^n\to \R^k$.  The derivative $Dh$ of $h$ is the restriction of the derivative $D\tilde{h}$ to $\R^n_u$.

\begin{definition}[{\cite[Definition~1.1]{alr}}]
Let $X_K$ be a second countable Hausdorff space and let $U$ be an open set in $X_K$.  An \textit{$n$-dimensional orbifold chart} over $U$ is a triple $(\tilde{U},\Gamma_U,\pi_U)$ satisfying
\begin{enumerate}
\item $\tilde{U}$ is a connected open set of $\R^n_{u}$ for some nonzero $u\in \R^n$,
\item $\Gamma_{U}$ is a finite group acting by diffeomorphisms on $\tilde{U}$,
\item $\pi_{U}\co\tilde{U}\to U$ is a continuous map such that $\pi_{U}\gamma=\pi_{U}$ for all $\gamma\in\Gamma_{U}$ and which induces a homeomorphism from $\tilde{U}/\Gamma_{U}$ onto $U$.
\end{enumerate}
Given two orbifold charts $(\tilde{U},\Gamma_U,\pi_U)$ and $(\tilde{V},\Gamma_V,\pi_V)$ with $U\subseteq V\subseteq X_K$, a (smooth) \textit{embedding} of orbifold charts $\kappa\co(\tilde{U},\Gamma_U,\pi_U)\hookrightarrow (\tilde{V},\Gamma_V,\pi_V)$ is a smooth embedding $\kappa\co\tilde{U}\hookrightarrow\tilde{V}$ with $\pi_V\kappa=\pi_U$.
\end{definition}

We now give the definition of a (smooth) orbifold with boundary. Note that all orbifolds considered in this paper will be smooth.

\begin{definition}
An $n$-dimensional (smooth) \textit{orbifold with boundary} $K$  is a second countable Hausdorff space $X_{K}$, called the underlying space of $K$, covered by a maximal atlas of $n$-dimensional orbifold charts $\{(U_a,\Gamma_{U_a},\pi_{U_a})\}_{a\in A}$ satisfying the following compatibility condition:  for any two charts $(\tilde{U}_{a}, \Gamma_{U_{a}},\pi_{U_{a}})$ and $(\tilde{U}_{b},\Gamma_{U_{b}},\pi_{U_{b}})$ and point $x\in U_{a}\cap U_{b}$ there is an open neighborhood $U_{c}\subset U_{a}\cap U_{b}$ about $x$ and a chart $(\tilde{U}_{c},\Gamma_{U_{c}}, \pi_{U_{c}})$ over $U_{c}$ such that there are smooth embeddings $\kappa_a\co (\tilde{U}_{c},\Gamma_{U_{c}},\pi_{U_{c}})\hookrightarrow (\tilde{U}_{a},\Gamma_{U_{a}},\pi_{U_{a}})$ and $\kappa_b\co(\tilde{U}_{c},\Gamma_{U_{c}},\pi_{U_{c}})\hookrightarrow(\tilde{U}_{b},\Gamma_{U_{b}},\pi_{U_{b}})$.  The \textit{boundary} $\partial K$ of $K$ is the set of all points $p\in X_K$ for which there is a chart $(\tilde{U}_a,\Gamma_{U_a},\pi_{U_a})$ such that $\langle \tilde{p}, u_a\rangle = 0$ for all points $\tilde{p}\in \pi_{U_a}^{-1}(p)$ along with the (smooth) orbifold structure inherited from $K$. Usually, it will convenient to abuse notation and refer to $X_K$ as $K$.
\end{definition}

The compatibility condition for orbifold charts allows us to consider the transition functions between overlapping charts $(\tilde{U}_{a},\Gamma_{U_{a}},\pi_{U_{a}})$ and $(\tilde{U}_{b},\Gamma_{U_{b}},\pi_{U_{b}})$, namely $\kappa_b\kappa_a^{-1}$ (see \cite[p.10]{alr}).  With transition functions in hand, we can define orientability.

\begin{definition}
We say that an orbifold $K$ is \textit{oriented} if it is covered by an atlas of orbifold charts such that for all pairs $(\tilde{U}_{a},\Gamma_{U_{a}},\pi_{U_{a}})$ and $(\tilde{U}_{b},\Gamma_{U_{b}},\pi_{U_{b}})$, $\det D(\kappa_b\kappa_a^{-1})>0$ where $D(\kappa_b\kappa_a^{-1})$ is the derivative of the transition function.
\end{definition}

Given an orbifold $K$, with or without boundary, we can construct a Riemannian structure on $K$ by placing a $\Gamma_U$-invariant Riemannian metric on the local cover $\tilde{U}$ of each orbifold chart $(\tilde{U},\Gamma_U,\pi_U)$ and patching these local metrics together with a partition of unity.  A smooth orbifold with a Riemannian structure is called a \textit{Riemannian orbifold} (see \cite{chiang,farsi}).

\begin{definition}
Let $(K,g)$ be an oriented Riemannian orbifold with boundary $\partial K$. Let $\nu$ be the outward pointing unit normal vector field along $\partial K$. A \emph{collar} is a smooth diffeomorphism $\Sigma\co \partial K\times [0,1)\to K$ onto an open neighborhood of $\partial K$ in $K$ such that $\Sigma(p,0)=p$ and the derivative $D\Sigma_{(p,0)}(0,1) = -\nu_p$ for all $p\in\partial K$. If we have fixed a collar $\Sigma$ on $K$, we say that $K$ is \textit{collared}.
\end{definition}

If $K$ is a compact oriented Riemannian orbifold that does not have singular points on the boundary, then $K$ admits a collar by \cite[Theorem~1.1.7]{schwarz}. The hypotheses that $K$ is oriented and Riemannian are required for $\nu$ to be defined. For the remainder of the paper, we assume that all Riemannian orbifolds are compact and oriented.  If the orbifold has boundary, we will assume that there are no singular points on the boundary and that the orbifold is collared.


\subsection{Connected sum metrics}
\label{subsec:ConnSumOflds}

Under certain conditions, it is possible to glue two Riemannian orbifolds with boundary together to obtain a connected sum orbifold.  Note that in this case, the metric on the connected sum will be piecewise differentiable, but not differentiable along the glued boundaries. By an \textit{orbifold domain} of an $n$-dimensional orbifold $O$, we mean a closed $n$-dimensional suborbifold with boundary. In this paper, we are primarily interested in the case of domains given by the closure of an open neighborhood of a singular point.

\begin{definition}
\label{def:ConnSum}
Let $O$ be a closed oriented orbifold, let $O_1$ be a closed orbifold domain with smooth boundary $\partial O_1$ such that there are no singular points on
$\partial O_1$. Let $O_2 = \overline{O\smallsetminus O_1}$
so that $\partial O_2 = \partial O_1$. We give $O_1$ and $O_2$ the orientations they inherit as domains in $O$. Let $g$ be a continuous metric on $O$ such that $g_1:=g_{|O_1}$ and
$g_2:=g_{|O_2}$ are both smooth. We refer to such a piecewise differentiable metric as a \emph{connected sum metric} on $O$.
\end{definition}

We note that throughout, whenever we glue two orbifolds with boundary together, we choose orientations on the two orbifolds in such a way that the connected sum has a natural orientation. We will define Sobolev spaces for connected sums in Section~\ref{subsec:SobolevSpacesConnectedSum}, define the Laplacian with respect to connected sum metrics in Section~\ref{subsec:LaplacianConnectedSum}, and use these definitions to study the spectrum of the Laplacian on the connected sum of two orbifolds $(O_1,g_1)$ and $(O_2,g_2)$ in Section~\ref{sec:ConnectedSums}.


\section{Sobolev spaces on orbifolds}
\label{sec:SobolevSpaces}

We now describe Sobolev spaces on orbifolds with boundary and orbifolds with connected sum metrics.  These Sobolev spaces will ultimately allow us to define for both settings the Hodge Laplacian.  In both cases, the Hodge Laplacian will be a non-negative, self-adjoint operator whose spectrum is discrete and tends to infinity.


\subsection{Preliminaries}
\label{subsec:SobolevSpacesPrelim}

As pointed out by Chiang \cite[p.320]{chiang}, given a smooth function $f$ on an oriented Riemannian orbifold $(K,g)$, with or without boundary, the lift of $f$ to any orbifold chart $\tilde{U}$ is $\Gamma_{U}$-invariant, and thus the gradient $\nabla_g f$ is $\Gamma_U$-invariant on $\tilde{U}$ as well.  Therefore, $\nabla_g f$ is well-defined on $(K,g)$. In the same way, the definitions of the exterior derivative $d$, exterior and interior products, the Hodge star operator $\star$, and the co-differential operator $\delta$ extend to the space $\Omega^p(K,g)$ of smooth differential $p$-forms on $K$.  We furthermore define the Laplacian $\Delta_g:\Omega^p(K,g)\to \Omega^p(K,g)$ by
\[
\Delta_g \omega =(d\delta+\delta d)\omega.
\]

We make sense of integration on orbifolds as follows.

\begin{definition}
Let $(K,g)$ be a compact oriented Riemannian orbifold, with or without boundary.  Suppose that
$\{(\tilde{U}_i, \Gamma_{U_i}, \pi_{U_i})\}^m_{i=1}$ is a finite covering of $(K,g)$ by orbifold charts.
Let $\{\rho_i\}_{i=1}^m$ be a partition of unity subordinate to $\{U_i\}^m_{i=1}$.
For a function $f\in \mathcal{C}^{\infty}(K)$ define the \textit{integral of $f$ over $K$} by
\[
\int_K f\, dv_g := \sum_{i=1}^m \frac{1}{\vert \Gamma_{U_i}\vert}\int_{\tilde{U}_i}\tilde{\rho}_i(\tilde{x})\tilde{f}(\tilde{x})\, d\tilde{v}_{\tilde{g}},
\]
where $\sim$ denotes the lift in all cases and $dv_g$ denotes the volume form with respect to $g$.  We note that although the definition makes use of a particular covering of $K$, it can be shown to be independent of the choice of covering (see \cite{chiang}, \cite{farsi} and \cite[p.34]{alr}).
Integrals of differential $p$-forms on $K$ are defined in the same way; see \cite[Section 2]{ChenRuan}.
\end{definition}


\subsection{Sobolev spaces on orbifolds with boundary}
\label{subsec:SobolevSpacesOfldsBdry}

Sobolev spaces of functions for closed orbifolds were introduced in \cite{chiang} and \cite{farsi}.  Although both authors gave the definition for orbifolds without boundary, the definition generalizes readily to differential forms on compact orbifolds with boundary.  (See \cite{schwarz} for a careful discussion of Sobolev spaces for manifolds with boundary.)

\begin{definition}
\label{def:Sobolev}
Let $(K,g)$ be an $n$-dimensional compact oriented Riemannian orbifold with boundary, and let $\Lambda^p(K,g)$ be the bundle of antisymmetric $p$-linear functions on the tangent space of $K$ with usual pointwise metric (see \cite[p.20]{schwarz}).  Note that $\Omega^p(K,g)$ is the space of smooth sections of $\Lambda^p(K,g)$.

Given $\omega\in \Omega^p(K,g)$, define $\vert \omega \vert_{J^0(\Lambda)}\in C^{\infty}(K)$ by
\[
\vert \omega\vert_{J^0(\Lambda)}:=\big(\langle \omega,\omega\rangle_{\Lambda^p(K,g)}\big)^{\frac{1}{2}}.
\]
Now, suppose $s$ is a positive integer and suppose that $E_i$ is a $g$-orthonormal frame on $(K,g)$.  Define $\vert \omega\vert_{J^s(\Lambda)}\in C^{\infty}(K)$ by
\[
\vert \omega\vert_{J^s(\Lambda)}:=\bigg(\vert\omega\vert_{J^{s-1}(\Lambda)}^2 + \sum\limits_{j=1}^n \vert\nabla_{E_j} \omega\vert_{J^{s-1}(\Lambda)}^2\bigg)^{\frac{1}{2}} \in\mathcal{C}^\infty(K),
\]
where $\nabla$ is the covariant derivative on forms induced by the Levi-Civita connection associated to $g$.  Set
\[
\Vert \omega\Vert_{W^{s,\ell}\Omega^p(K,g)} =  \bigg(\int_K \vert \omega\vert_{J^s(\Lambda)}^{\ell}\, dv_g\bigg)^{\frac{1}{\ell}}.
\]
The \textit{Sobolev space} $W^{s,\ell}\Omega^p(K,g)$  is the completion of $\Omega^p(K,g)$ with respect to the norm $\Vert\cdot\Vert_{W^{s,\ell}\Omega^p(K,g)}$.  For the case $s=0$, we denote $W^{0,\ell}\Omega^p(K,g)$ by $L^{\ell}\Omega^p(K,g)$ and for $\ell=2$,  we denote $W^{s,2}\Omega^p(K,g)$ by $H^s\Omega^p(K,g)$.
When $p = 0$, we simply denote $H^s(K,g) = H^s\Omega^0(K,g)$ and $L^{\ell}(K,g) = L^{\ell}\Omega^0(K,g)$.
\end{definition}

We remark that $H^0\Omega^p(K,g) = L^2\Omega^p(K,g)$ and that the $H^1$-norm is given on functions by
\[
\Vert f\Vert_{H^1(K,g)} =\Big( \langle f,f\rangle_{L^2(K,g)} +\langle \nabla_g f,\nabla_g f\rangle_{L^2(K,g)}\Big)^{\frac{1}{2}},
\]
where $\nabla_g$ denotes the gradient.
Furthermore, by taking the $W^{s,\ell}$-limit of a sequence of smooth forms, we can extend the notions listed at the beginning of Section~\ref{subsec:SobolevSpacesPrelim}, e.g., the exterior derivative $d$, integration, and the tangential and normal components of a form, to Sobolev spaces having the appropriate level of regularity for the particular notion under consideration. In addition, the definition of $W^{s,\ell}\Omega^p(K,g)$ and hence $H^s\Omega^p(K,g)$ can be extended to non-integer values of $s$ in the standard way using the Fourier transform; see \cite[Section~4.1]{taylor}.

As usual, $H^s\Omega^p(K,g) \subset H^t\Omega^p(K,g)$ when $s>t$ and the following orbifold version of Rellich's theorem holds. The proof of the same result for manifolds \cite[Theorem 1.3.6]{schwarz} is based on local arguments, and hence carries directly to the orbifold setting by restricting to locally defined forms that are invariant under the action of a finite group. See also \cite[Theorem 2.1]{chiang} or \cite[Theorem 2.4]{farsi} for the case $p=0$.

\begin{theorem}[Rellich's theorem]
\label{thrm:rellich}
Let $(K,g)$ be a compact oriented Riemannian orbifold. If $(K,g)$ has boundary $\partial K$, assume that $\partial K$ is a manifold, i.e., contains no singular points.  Suppose that $s>t$.  Then the inclusion $H^s\Omega^p(K,g)$ into $H^t\Omega^p(K,g)$ is compact.
\end{theorem}

The trace theorem for manifolds with boundary is proven by patching together local trace operators between Sobolev spaces on regions in $\R^n$ and their boundaries, and therefore carries over directly to the case of collared orbifolds with no singular points on the boundary. See \cite[Theorem 1.3.7]{schwarz}; see also \cite[p.258]{evans} and \cite[Chapter 4, Propositions~1.6 and 4.5]{taylor} for the case $p=0$. Note that when $p > 0$, the restriction $\omega_{|\partial K}$ of a differential $p$-form $\omega$ on $K$ to the boundary is a section of the bundle $\Lambda^p(K,g)_{|\partial K}$ and not a differential form on $\partial K$. The Sobolev norm $\Vert\omega_{|\partial K}\Vert_{H^s\Omega^p(K,g)_{|\partial K}}$ of such a section is defined as in Definition~\ref{def:Sobolev} in terms of the metric that $\Lambda^p(K,g)_{|\partial K}$ inherits from the metric on $\Lambda^p(K,g)$ and by integrating over $\partial K$, see \cite[Definition 1.3.1]{schwarz}. Of course, when $p = 0$, $\omega$ restricts to a function on $\partial K$ with the usual Sobolev norm.

\begin{theorem}[Trace theorem]
\label{thrm:trace}
Let $(K,g)$ be a compact oriented Riemannian orbifold with nonempty manifold boundary $\partial K$.
\begin{enumerate}
\item   For $s=0,1$, the trace map $\omega\mapsto\omega_{\vert\partial K}$ is a compact continuous operator $H^{s+1}\Omega^p(K,g)$ to $H^s\Omega^p(K,g)_{|\partial K}$, i.e., $\Vert \omega_{\vert\partial K}\Vert_{H^s\Omega^p(K,g)_{|\partial K}}\leq C \Vert \omega\Vert_{H^{s+1}\Omega^p(K,g)}$ for some choice of $C$ independent of $\omega$.
\item   The trace map extends uniquely to a continuous operator $H^{1}\Omega^p(K,g)\to H^{1/2}\Omega^p(K,g)_{|\partial K}$.
\end{enumerate}
\end{theorem}

For an orbifold $(K,g)$ with boundary, let $\partial g$ denote the induced metric on $\partial K$.
Given a vector field $X$ on $K$, considering the restriction of $X$ to the boundary $\partial K$, we can decompose $X=X^{\Vert}+X^{\perp}$ into its tangential and normal components, respectively. As in \cite{schwarz}, we then define
\[
\bft\omega(X_1,\dots, X_p)=\omega(X_1^{\Vert},\dots,X_p^{\Vert})\hspace{8mm} \textrm{and}\hspace{8mm} \bfn\omega = \omega_{\vert\partial K}-\bft\omega.
\]
The manifold version of Green's theorem is derived from Stokes' theorem by arguments that are all local in nature, and thus generalize to orbifolds; see, for example, \cite[Chapter 2, Proposition 1.2]{schwarz} and \cite[Proposition 2.4.1]{taylor}. As the proof of Stokes' theorem is also given locally and patched together with a partition of unity, see \cite[Chapter 1, Proposition 13.3]{taylor}, the proof of Green's theorem extends directly to the orbifold case, yielding the following.

\begin{theorem}[Green's theorem]
\label{thrm:Green}
On a compact oriented Riemannian orbifold $(K,g)$ with manifold boundary $\partial K$, let $\omega\in H^1\Omega^{p-1}(K,g)$ and $\eta\in H^1\Omega^p(K,g)$. Then
\begin{enumerate}
\item[(1)]
$\displaystyle\langle d\omega, \eta\rangle_{L^2\Omega^p(K,g)} =\langle\omega,\delta\eta\rangle_{L^2\Omega^{p-1}(K,g)}
+ \int_{\partial K} \bft\omega\wedge\star \bfn\eta$.
\end{enumerate}
For $f, h\in H^2(K,g)$, we have
\smallskip
\begin{enumerate}
\item[(2)]
$\displaystyle\langle \Delta_g f, h\rangle_{L^2(K,g)} =\langle \nabla_g f,\nabla_g h\rangle_{L^2(K,g)}
+ \int_{\partial K}h\big(\nu(f)\big)\, dv_{\partial g}$.
\smallskip
\item[(3)]
$\displaystyle\langle \Delta_g f, h\rangle_{L^2(K,g)} -\langle f,\Delta_g h\rangle_{L^2(K,g)}= \int_{\partial K}\Big(h\big(\nu(f)\big)-f\big(\nu(h)\big)\Big)\, dv_{\partial g}$,
\end{enumerate}
where $\nu$ denotes the unit outward pointing normal on $\partial K$.
\end{theorem}
We note that there are varying conventions regarding the sign of the Laplace operator and the choice of inward or outward normal, leading to formulations of Green's theorem that may differ by a minus sign on individual terms from the one we have presented here.


\subsection{Sobolev spaces on orbifolds with connected sum metrics}
\label{subsec:SobolevSpacesConnectedSum}

Let $O$ be a closed oriented orbifold with connected sum metric as in Definition~\ref{def:ConnSum}.  By Theorem~\ref{thrm:trace}~(1), the following definitions make sense.   See \cite[Definitions 1.1 and 1.4]{annecolbois} and \cite[page 17]{takahashiGapPForms}.

\begin{definition}
\label{def:ConnSumSobolev}
We define $L^2\Omega^p(O,g)$ by
\[
    L^2\Omega^p(O, g) := L^2\Omega^p(O_1,g_1)\oplus L^2\Omega^p(O_2,g_2)
\]
with the componentwise inner product. Note that in this way, $L^2\Omega^p(O,g)$ is a Hilbert space.

The \textit{first Sobolev space} $H^1\Omega^p(O,g)$ is the subspace of $H^1\Omega^p(O_1,g_1)\oplus H^1\Omega^p(O_2,g_2)$ defined by
\begin{align*}
H^1\Omega^p(O,g):=\Big\{\omega &=(\omega_1,\omega_2)\in H^1\Omega^p(O_1,g_1)\oplus H^1\Omega^p(O_2,g_2)\,\vert\\
&\bft_1\omega_{1\vert\partial O_1}=\bft_2\omega_{2\vert\partial O_2}\text{ in } L^2\Omega^p(\partial O_1,\partial g_1),\\ &\star\bfn_1\omega_1=-\star\bfn_2\omega_2 \text{ in } L^2\Omega^{n-p}(\partial O_1,\partial g_1),\\
&d\omega\in L^2\Omega^{p+1}(O,g), \text{ and } \delta\omega \in L^2\Omega^{p-1}(O,g)\Big\}.
\end{align*}
The inner product on $H^1\Omega^p(O,g)$ is given by the direct sum of inner products on the two component spaces.

The \textit{second Sobolev space} $H^2\Omega^p(O,g)$ is the subspace of $H^2\Omega^p(O_1,g_1)\oplus H^2\Omega^p(O_2,g_2)$ defined by
\begin{align*}
H^2\Omega^p(O,g) :=\Big\{\omega &=(\omega_1,\omega_2)\in H^2\Omega^p(O_1,g_1)\oplus H^2\Omega^p(O_2,g_2)\,\vert\\
&\bft_1\omega_{1\vert\partial O_1}=\bft_2\omega_{2\vert\partial O_2}\text{ in } H^1\Omega^p(\partial O_1,\partial g_1),\\
&\star\bfn_1\omega_1=-\star\bfn_2\omega_2\text{ in } H^1\Omega^{n-p}(\partial O_1,\partial g_1),\\
&d\omega\in H^1\Omega^{p+1}(O,g), \text{ and } \delta\omega \in H^1\Omega^{p-1}(O,g)\Big\}.
\end{align*}
where $d\omega=(d\omega_1, d\omega_2)$ and $\delta\omega=(\delta\omega_1,\delta\omega_2)$.
The inner product on $H^2\Omega^p(O,g)$ is given by the direct sum of the inner products on the two component spaces.
\end{definition}

We note that for $p=0$ and $s = 1$ or $2$, Definition~\ref{def:ConnSumSobolev} yields
\begin{align*}
\label{eq:Sobolev-functions-1}
    H^s(O,g) := \Big\{f &=(f_1,f_2)\in H^1(O_1,g_1)\oplus H^1(O_2,g_2)\, \vert\,\\
    &f_{1\vert\partial O_1} = f_{2\vert\partial O_2} \text{ in }L^2(\partial O_1,\partial g_1),\\
    &(\nu_1(f_1))_{\vert\partial O_1}
        =-(\nu_2(f_2))_{\vert\partial O_2}\text{ in }L^2(\partial O_1, \partial g_1).\\
    &\text{and }df\in H^{s-1}(O,g) \Big\},
\end{align*}
where $\nu_1$ and $\nu_2$ denote the outward unit normal vector fields along $\partial O_1$ and
$\partial O_2$, respectively.

We have the following version of Rellich's theorem for orbifolds that carry connected sum metrics.

\begin{theorem}[Rellich's theorem for connected sum metrics]
\label{thrm:RellichPcSm}
For $L^2\Omega^p(O,g)$ and $H^1\Omega^p(O,g)$ as above, the inclusion of $H^1\Omega^p(O,g)$ into $L^2\Omega^p(O,g)$ is compact.
\end{theorem}
\begin{proof}
By assumption, the orbifolds $(O_1, g_1)$ and $(O_2,g_2)$ are collared orbifolds with no singular points on the boundary.  Thus, by Theorem~\ref{thrm:rellich}, $H^1\Omega^p(O_1,g_1)$ is compactly embedded in $L^2\Omega^p(O_1,g_1)$ and $H^1\Omega^p(O_2,g_2)$ is compactly embedded in $L^2(O_2,g_2)$.  But the norms on these respective spaces are all computed via a direct sum, and since $H^1\Omega^p(O_1,g_1)$ is a closed subspace of $H^1\Omega^p(O_1,g_1)\oplus H^1\Omega^p(O_2,g_2)$, it follows directly that the inclusion of $H^1\Omega^p(O,g)$ into $L^2\Omega^p(O,g)$ is compact.
\end{proof}


\subsection{Quadratic forms}
\label{subsec:Quad}

Recall that a \emph{quadratic form} is determined by a triple $(\mathcal{H}, D(q), q)$ where
$\mathcal{H}$ is a Hilbert space; $D(q)\subseteq \mathcal{H}$ is a closed subspace of $\mathcal{H}$,
the domain of $q$; and $q(x) = f(x,x)$ for a symmetric non-negative bilinear form  $f: D(q)\times D(q)\to\R$.
The norm of the  quadratic form $q$ is defined to be
\begin{equation}
\label{eq:QuadFormNorm}
    \Vert q\Vert=	\sup_{\{x \in D(q)\,\vert\, \langle x,x\rangle=1\}} q(x),
\end{equation}
where $\langle\cdot,\cdot\rangle$ denotes the inner product on $\mathcal{H}$.

We note that we can define the spectrum of a quadratic form using Rayleigh-Ritz variational formulae.  Further, when the quadratic form induces a non-negative self-adjoint operator (e.g., the Laplacian), the spectrum of the quadratic form coincides with that of the operator, see \cite[Chapter 4]{davies}.

We now introduce a bilinear form $q$ on $H^1\Omega^p(K,g)$ and corresponding quadratic form that induces the Laplacian on $(K,g)$. This form will be instrumental in the proof of the Hodge decomposition theorem, Theorem~\ref{thrm:Hodge}, as well as several arguments thereafter. It is an essential ingredient to the proofs in Section~\ref{sec:Laplacian} that the spectrum of the Laplacian in our context has the usual properties of being discrete and tending to infinity---see in particular Section~\ref{subsec:LaplacianBoundary}---as well as in Sections~\ref{sec:SpecDefDim3Func} and \ref{sec:SpecDefDim3Form}, where it allows us to study the spectrum via Rayleigh-Ritz variational formulae.

Suppose that $\omega,\eta\in H^1\Omega^p(K,g)$.  Let
\begin{equation}
\label{eqn:qform}
    q(\omega,\eta)= \langle d\omega,d\eta\rangle_{L^2\Omega^{p+1}(K,g)}+\langle\delta\omega,\delta\eta\rangle_{L^2\Omega^{p-1}(K,g)},
\end{equation}
and note that
\[
    q(\omega,\omega)=\Vert d\omega\Vert_{L^2\Omega^{p+1}(K,g)}+\Vert \delta\omega\Vert_{L^2\Omega^{p-1}(K,g)}.
\]
In a slight abuse of notation, we use $q$ to denote the bilinear form defined in Equation~\eqref{eqn:qform} as well as the quadratic form $q(\omega)=q(\omega,\omega)$ defined on $D(q)=H^1\Omega^p(K,g)$ in $\mathcal{H}=L^2\Omega^p(K,g)$.
By Green's Theorem (Theorem~\ref{thrm:Green}), we immediately have the following (see \cite[Corollary 2.1.4]{schwarz}).

\begin{proposition}
\label{prop:qOrbBdry}
Let $(K,g)$ be a compact oriented Riemannian orbifold with manifold boundary $\partial K$.
For $\omega \in H^2\Omega^p(K, g)$ and $\eta\in H^1\Omega^p(K, g)$,
the bilinear form $q$ is given by
\[
    q(\omega,\eta)  = \langle \Delta_g \omega, \eta\rangle_{L^2\Omega^p(K,g)}
        + \int\limits_{\partial K} \bft \eta \wedge \star\bfn d\omega
        - \int\limits_{\partial K} \bft \delta \omega \wedge \star\bfn \eta.
\]
\end{proposition}


\section{The Hodge decomposition theorem for orbifolds with boundary}
\label{sec:Hodge}

In this section, we state and prove the Hodge decomposition theorem (Theorem~\ref{thrm:Hodge}) for compact orbifolds with manifold boundary. Our exposition closely follows that of Schwarz \cite[Chapters 1,2]{schwarz} and his proof of the corresponding result for manifolds with boundary, \cite[Theorem 2.4.2]{schwarz}; see also \cite{bailey} for the case of orbifolds without boundary, and recall that orbifolds are also called ``$V$-manifolds." Indeed, as we will outline, practically all the proofs in \cite[Chapters 1,2]{schwarz} extend  verbatim to the orbifold case.
However, since these extensions are not noted in the literature, we will give a brief account of them
so that they will be available for our and future use.

Although we expect that a  similar  Hodge decomposition theorem is likely  to hold for compact orbifolds with singular points on the boundary, additional technical
difficulties arise, especially related to the extension of the results of \cite[Chapter XX]{hormander}.
More precisely, the proof that the differential operator
associated to an elliptic boundary value problem
is Fredholm, which is used for proving the completeness of the Hodge decomposition,
uses the Douglis-Nirenberg theory on the boundary,
see \cite[Theorem 20.1.3, Corollary 20.1.4, Proposition 20.1.5]{hormander}. Since such theory
has not explicitly been developed for orbifolds, this makes the extension of the
Hodge decomposition to the case
with singularities on the boundary more involved.
Therefore, since our connected sum constructions are always applied to the case of a
manifold boundary, we have restricted our attention to this case.

We start with some basic definitions.

\begin{definition}
\label{def:BoundaryConditions}
Let $(K,g)$ be an orbifold with manifold boundary, and suppose that $\omega\in H^1\Omega^p(K,g)$.   We say that $\omega$ satisfies \textit{Dirichlet boundary conditions} if $\bft\omega=0$ on $\partial K$, and we say that $\omega$ satisfies \textit{Neumann boundary conditions} if $\bfn\omega=0$ on $\partial K$.

For $\omega\in H^2\Omega^p(K,g)$, we say that $\omega$ satisfies \textit{relative} boundary conditions if $\bft \omega=0$ and $\bft\delta\omega=0$.  We say that $\omega$ satisfies \textit{absolute} boundary conditions if $\bfn\omega=0$ and $\bfn d\omega=0$ (see~\cite[p.728]{mcgowan} and cf.~Equations~\eqref{eq:boundary-value-problem}). We note that for the case of functions, relative boundary conditions are equivalent to Dirichlet boundary conditions and absolute boundary conditions are equivalent to Neumann boundary conditions.
\end{definition}

Similarly, we define the following subspaces of $L^2\Omega^p(K,g)$.

\begin{definition}
Let $H^1\Omega_D^p(K,g)$ and $H^1\Omega_N^p(K,g)$ be the subspaces of $H^1\Omega^p(K,g)$ with vanishing tangential and normal components respectively.  Let
\[
\mathcal{E}^p(K)=\big\{d\alpha\,\vert\, \alpha\in H^1\Omega_D^{p-1}(K,g)\big\}
\]
be the space of exact $p$-forms with Dirichlet conditions, and let
\[
\mathcal{C}^p(K)=\big\{\delta\beta\,\vert\, \beta\in H^1\Omega_N^{p+1}(K,g)\big\}
\]
be the space of coexact $p$-forms with Neumann conditions.  Set $\mathcal{E}^0(K)=\{0\}$ and $\mathcal{C}^d(K) =\{0\}$.
Finally, let $\mathcal{H}^p(K)$ denote the space of \textit{harmonic fields} in $H^1\Omega^p(K,g)$, i.e., the set
\[
\mathcal{H}^p(K)=\big\{\omega\in H^1\Omega^p(K,g)\,\vert\, d\omega=0\text{ and } \delta\omega=0\big\},
\]
and let $L^2\mathcal{H}^p(K)$ be the $L^2$-closure of $\mathcal{H}^p(K)$.
\end{definition}

\begin{remark}
We note that we use the term \textit{harmonic fields} to denote elements $\omega\in H^1\Omega^p(K,g)$ for which $d\omega=\delta\omega=0$ and we use the term \textit{harmonic forms} to denote forms for which $\Delta_g\omega=0$.
\end{remark}

We now state our main theorem.  See \cite[Theorem 2.4.2]{schwarz} for the manifold version. The proof of this theorem occupies the rest of this section.

\begin{theorem}[Hodge decomposition theorem for orbifolds with boundary]
\label{thrm:Hodge}
Let $(K,g)$ be a compact oriented Riemannian orbifold with manifold boundary.   Then $L^2\Omega^p(K,g)$ splits into the $L^2$-orthogonal direct sum
\[
    L^2\Omega^p(K,g)=\mathcal{E}^p(K)\oplus \mathcal{C}^p(K)\oplus L^2\mathcal{H}^p(K).
\]
\end{theorem}

The proof of Theorem~\ref{thrm:Hodge} depends on a number of preliminary results which we present below. By combining these results, the proof of Theorem~\ref{thrm:Hodge} will follow naturally at the end of this section.

We begin by noting that the proof of Gaffney's inequality extends to the orbifold case which, along with
the definition of the bilinear form $q$ in Equation~\eqref{eqn:qform} and Proposition~\ref{prop:qOrbBdry}, yields the following.

\begin{lemma}[{Gaffney's inequality \cite[Corollary 2.1.6]{schwarz}}]
\label{lem:Gaffney}
Let $(K,g)$ be a compact oriented Riemannian orbifold with manifold boundary. Then there is a constant $C$ depending only on $(K,g)$ such that for each $\omega \in H^1\Omega^p(K, g)$ with $\bft\omega = 0$, we have
\[
    \Vert\omega\Vert_{H^1\Omega^p(K, g)}^2 \leq C \big( q(\omega,\omega)  + \Vert\omega\Vert_{L^2\Omega^p(K, g)}^2\big).
\]
\end{lemma}

\begin{definition}
Let $\HH_D^p(K)=\mathcal{H}^p(K)\cap H^1\Omega_D^p(K,g)$, the subspace of harmonic fields in $H^1 \Omega_D^p(K,g)$, and let
$\HHDkK^{\obot}$ denote the $L^2$-orthogonal complement of $\HHDkK$ in $H^1 \Omega_D^p(K,g)$.
Define subspaces $\mathcal{H}_N^p(K)$ and $\mathcal{H}_N^p(K)^{\obot}$ of $H^1\Omega_N^p(K,g)$ analogously.
We refer to elements of $\HH_D^p(K)$ as \emph{Dirichlet fields} and elements of $\mathcal{H}_N^p(K)$ as \emph{Neumann fields}.
\end{definition}

We make note as well of the following relationship between $q$ and the $H^1$-norm on a subspace of $H^1\Omega^p(K,g)$.
Using Lemma~\ref{lem:Gaffney} and Theorem~\ref{thrm:rellich}, the proof of \cite[Proposition 2.2.3]{schwarz} generalizes directly to orbifolds.

\begin{proposition}[{\cite[Proposition 2.2.3]{schwarz}}]
\label{prop:qH1norms}
Suppose that $\omega\in \HHDkKo$.  There exist positive constants $c,C$ such that
\[
c\Vert\omega\Vert_{H^1\Omega^p(K,g)}^2 \leq q(\omega,\omega)\leq C\Vert \omega\Vert_{H^1\Omega^p(K,g)}^2.
\]
A similar relationship holds for $\omega\in \HHNkKo$.
\end{proposition}


\subsection{The Dirichlet and Neumann potentials and their regularity}
\label{sub:Dirichlet-pot-regularity}

We now establish technical results involving the existence and regularity of Dirichlet and Neumann potentials for elements of $\HHDkK^{\perp}$ and $\mathcal{H}_N^p(K)^{\perp}$ respectively.  The proof of the Hodge decomposition theorem for orbifolds with boundary (Theorem~\ref{thrm:Hodge}) will be presented at the end of Section~\ref{sec:Hodge}, following these results.


\subsubsection{Definition of the Dirichlet and Neumann potentials}

Using Proposition~\ref{prop:qH1norms}, the proof of the existence of the Dirichlet and Neumann potentials follows in a manner identical to \cite[Theorem 2.2.4 and p.72]{schwarz}.

\begin{theorem}[{\cite[Theorem 2.2.4]{schwarz}}]
\label{def:dirichlet-potential}
Let $(K, g)$  be a compact oriented Riemannian orbifold with manifold boundary.
For  each $\eta \in \HHDkK^{\perp}, $ there exists a differential
form $\DP \in \HHDkK^{\obot}$ such that
\begin{equation}
q(\DP, \xi)
= \langle \eta, \xi \rangle_{L^2\Omega^p(K,g)},\ \forall \xi \in  \HOKE .
\end{equation}
The differential form $\DP$ is uniquely determined by $\eta$ and is called the Dirichlet
potential of $\eta$. Conversely, $\eta$ is uniquely determined by its Dirichlet potential $\DP$.   A similar form $\varphi_N\in\mathcal{H}_N^p(K)^{\obot}$ exists for each $\eta\in \mathcal{H}_N^p(K)^{\perp}$, and is called the Neumann potential for $\eta$.
\end{theorem}

\begin{remark}
\label{rmk:potential-weak-solution}
Note that although we do not currently know that $\DP$ is sufficiently differentiable, we can formally apply Proposition~\ref{prop:qOrbBdry} to the result of Theorem~\ref{def:dirichlet-potential} to see that the Dirichlet potential $\DP$ weakly solves the boundary value problem
\begin{equation}
\label{eq:boundary-value-problem}
\begin{split}
    \Delta_g\DP = \eta   \quad &\text{on }  K,     \\
    \bft \DP = 0 \hbox{ and }  \bft \delta \DP = 0  \quad &\text{on }  \partial K ,
\end{split}
\end{equation}
i.e. $\Delta_g\DP = \eta$ with relative boundary conditions, see Definition~\ref{def:BoundaryConditions}.
\end{remark}

This remark is strengthened in \cite[Theorem 2.2.5]{schwarz} by showing that the Dirichlet potential is also a strong solution of Equation~\eqref{eq:boundary-value-problem}. We summarize the extension of this and other relevant results to orbifolds in Section~\ref{subsubsec:Regularity}. Because the argument for the regularity of $\DP$ we provide is very similar to the manifold case proved in \cite[Section 2.3]{schwarz}, we will only give a sketch while outlining the few differences. In particular,  we will need to use results on integration of vector fields on orbifolds, as well as standard (pseudo-)differential operators techniques.

Finally, again as in the comment \cite[p.72]{schwarz}, all results that we state and prove here for the Dirichlet potential hold for the Neumann potential $\varphi_N$ as well.  See \cite[Theorem 2.2.7]{schwarz} and Theorem~\ref{thrm:neumann-potential}.


\subsubsection{Regularity on the interior and on the boundary of the Dirichlet and Neumann potentials}
\label{subsubsec:Regularity}

The main new ingredient in the orbifold case are the results of \cite{hepworth} on flows of vector fields.
Let $X$ be a compactly supported vector field on the orbifold $K$ with boundary that is either supported away from the boundary or such that $X_{\vert \partial K}$ is parallel to $\PK$. Then $X$ can be integrated to produce a global flow $\FLOW$ using \cite[Theorem 5.13]{hepworth} applied to the double of $K$. Then we proceed as in \cite[Section 2.3]{schwarz} in the orbifold case.
Define the  operator
\begin{align}
\SIGMA: \LOkK \to \LOkK,\ \SIGMA &:=  \frac1t \Big( \PB  -\text{Id}         \Big),
\end{align}
and then the proof of \cite[Lemma 2.3.1]{schwarz} extends verbatim to the orbifold case.

We now outline the proof of ${H}^2$-regularity of the Dirichlet potential. We  will need to prove ${H}^2$-regularity separately at interior points and boundary points. The proof of \cite[Lemma 2.3.2]{schwarz},
which demonstrates ${H}^2$-regularity in the interior as well as some covariant derivative ${H}^1$-regularity on the boundary, is identical in the orbifold case using results in \cite{hepworth} on integration of orbifold vector fields as described above. It also involves \cite[Proof of Lemma 2.3.1]{schwarz} and   \cite[Equation (1.3.15)]{schwarz} (Meyers-Serrin's theorem, i.e., \cite[Theorem 1.3.6(a)]{schwarz} or \cite[Theorem 1.3.3(a)]{schwarz}), which are proven identically for orbifolds. The remaining ingredients hold in general for Hilbert spaces. Then \cite[Lemma 2.3.3]{schwarz},
which demonstrates ${H}^2$-regularity of the Dirichlet potential on the boundary, also follows as in the manifold case. The proof uses orbifold analogues of the Bochner-Weitzenb\"ock formula \cite[Equation (1.2.23)]{schwarz}, \cite[Lemma 2.3.2]{schwarz}, and the Meyers-Serrin theorem, \cite[Equation (1.3.15)]{schwarz}, \cite[Theorem 1.3.6(a)]{schwarz}, and \cite[Theorem 1.3.3(a)]{schwarz}, which extend to orbifolds using local charts. Then the proof of \cite[Theorem 2.2.5]{schwarz}, which strengthens Remark~\ref{rmk:potential-weak-solution}, applies directly to orbifolds.

In order to prove the completeness of the Hodge decomposition, we need \cite[Theorem 1.6.2]{schwarz} and Theorem~\ref{thrm:theorem2.2.6} below, which we have confirmed for orbifolds using the preceding results and standard results on  elliptic (pseudo-)differential operators.  See \cite[Chapter XX]{hormander} and \cite[Theorem 1.6.2, Theorem 2.2.6, Corollary 2.3.5, and p.78]{schwarz}. In particular we note that to prove \cite[Theorem 1.6.2]{schwarz} for orbifolds,
\cite[Theorem 20.1.2]{hormander}, \cite[Theorem 20.1.8]{hormander} and \cite[Proposition 20.1.11]{hormander} need to be extended to orbifolds. This can be done by using a variety of methods (depending on the context),  which include using  a collared neighborhood of the boundary of the form  $\partial K \times [0,1)$ (where $\partial K$ has no singular point),  using local orbifold charts, as well as the theory of pseudo-differential operators on compact orbifolds. For example, the existence of parametrices is shown in \cite{Bucicovschi}.
We then have the following for relative boundary conditions. For the argument, which carries directly to the orbifold setting, see \cite[Theorem 2.2.6 and Corollary 2.3.5]{schwarz}.

\begin{theorem}
\label{thrm:theorem2.2.6}
Let $(K,g)$ be a compact oriented Riemannian orbifold with manifold boundary.
\begin{enumerate}
\item If $\psi$ is a Dirichlet field, $\psi\in \mathcal{H}_D^p(K)$, then $\psi$ is smooth.
\item If $\eta\in \mathcal{H}_D^p(K)^{\perp}$ is of Sobolev class $W^{s,p}$ (with $s\in \mathbb{N}_0$ and $p\geq 2$), then the Dirichlet potential $\varphi_D$ is in $W^{s+2,p}\Omega^p(K)$.
\end{enumerate}
\end{theorem}

By the remarks \cite[p.72]{schwarz}, we also have the following for absolute boundary conditions.  See \cite[Theorem 2.2.7]{schwarz}.

\begin{theorem}
\label{thrm:neumann-potential}
Let $(K,g)$ be a compact oriented Riemannian orbifold with manifold boundary.
\begin{enumerate}
\item The space $\mathcal{H}^p_N(K)$ of Neumann fields is finite dimensional, and consists of smooth differential forms $\psi\in\Omega^p(K,g)$.
\item For each $\eta\in \mathcal{H}_N^p(K)^{\perp}$ there exists a unique Neumann potential $\varphi_N\in\mathcal{H}_N^p(K)^{\obot}\cap H^2\Omega^p(K,g)$.
\item Suppose that $s\in\mathbb{N}_0$ and $p\geq 2$.  If $\eta\in \mathcal{H}_N^k(K)^{\perp}\cap W^{s,p}\Omega^p(K)$, the Neumann potential $\varphi_N$ is of Sobolev class $W^{s+2,p}(K)$.
\end{enumerate}
\end{theorem}

The proof of Theorem~\ref{thrm:Hodge} now follows.

\begin{proof}[Proof of Theorem~\ref{thrm:Hodge}]
Since \cite[Theorems 1.6.2 and 2.2.5]{schwarz} hold for orbifolds as noted above, using Theorems~\ref{thrm:theorem2.2.6} and \ref{thrm:neumann-potential}, the proof of Theorem~\ref{thrm:Hodge} follows exactly as the proof of Theorem 2.4.2 in \cite{schwarz}.
\end{proof}

\begin{remark}
All of the results and proofs above hold if the boundary of $(K,g)$ is empty (and hence with vacuous boundary conditions).  Thus, the Hodge decomposition theorem as in \cite{bailey} holds for closed oriented orbifolds as well.
\end{remark}


\section{The Laplacian on orbifolds}
\label{sec:Laplacian}

In this section, we confirm that the spectrum of the Laplacian acting on $p$-forms has the familiar properties of being discrete and tending to infinity, both in the case of orbifolds with boundary under either Dirichlet or Neumann conditions, and in the case of orbifolds with connected sum metrics.


\subsection{The Laplacian on orbifolds with boundary}
\label{subsec:LaplacianBoundary}

Suppose that $(K,g)$ is an $n$-dimensional compact oriented Riemannian orbifold with manifold boundary.  As in Equation~\eqref{eqn:qform}, we define a bilinear form $q$ on $H^1\Omega^p(K,g)$ and recall that $q$ induces the Laplacian on $H^2\Omega^p(K,g)$, with  boundary terms as in Proposition~\ref{prop:qOrbBdry}.

Recall from Definition~\ref{def:BoundaryConditions} that $\omega\in H^2\Omega^p(K,g)$ satisfies \textit{relative} boundary conditions if $\bft \omega=0$ and $\bft\delta\omega=0$ and \textit{absolute} boundary conditions if $\bfn\omega=0$ and $\bfn d\omega=0$. We use the following notation for the nondecreasing sequences of nonzero eigenvalues including repetitions,
$k=1,2,\ldots$:
\begin{itemize}
\item   $\lambda_{p,k}^D(K,g)$, to denote the $k$th eigenvalue of the Laplacian acting on $p$-forms with relative boundary conditions,
\item   $\lambda_{p,k}^N(K,g)$, to denote the $k$th eigenvalue of the Laplacian acting on $p$-forms with absolute boundary conditions, and
\item   $\mu_{p,k}^N(K,g)$, to denote the $k$th eigenvalue of the Laplacian acting on coexact $p$-forms with absolute boundary conditions.
\end{itemize}
When $K$ has empty boundary, we use
\begin{itemize}
\item   $\lambda_{p,k}(K,g)$, to denote the $k$th eigenvalue of the Laplacian acting on $p$-forms and
\item   $\mu_{p,k}(K,g)$, to denote the $k$th eigenvalue of the Laplacian acting on coexact $p$-forms.
\end{itemize}
When $0$ occurs as an eigenvalue, it has index $k = 0$. See Section~\ref{sec:SpecDefDim3Form}.
We will omit $(K,g)$ when it is clear from the context. When restricting to the case of the Laplacian acting on functions, we will omit the subscript $p=0$, i.e., use $\lambda_k(K,g) = \lambda_{0,k}(K,g)$, etc.

\begin{remark}
\label{rem:ExactCoexactSpec}
Since for each $p=1,\dots, n$, the exterior derivative $d$ maps the space of coexact $(p-1)$-eigenforms with absolute boundary conditions isomorphically onto the space of exact $p$-eigenforms with absolute boundary conditions (see \cite[p.34]{gornet-mcgowan}), $\mu_{p,k}^N(K,g)$ is also the $k$th eigenvalue of the Laplacian acting on exact $(p+1)$-forms with absolute boundary conditions. Then by Theorem~\ref{thrm:Hodge}, the eigenvalue spectrum on $p$-forms is the union of the eigenvalue spectrum on the coexact $p$-forms and coexact $(p-1)$-forms, all with absolute boundary conditions. When $K$ is closed, we also have by Hodge duality that $\mu_{p,k}(K,g)=\mu_{n-p-1,k}(K,g)$. Therefore, in the case that $K$ has no boundary, by understanding the behavior of coexact $p$-eigenforms for $p=0,1,\dots, \lfloor (n-1)/2\rfloor$, we understand the behavior of the spectrum of the Laplacian on $p$-forms for all $p=0,1,\dots, n$, see \cite[p.968]{JammesPrecrip}.
\end{remark}

Suppose that $\omega,\eta\in H^2\Omega^p(K,g)$, with either relative or absolute boundary conditions.  Then $\langle \Delta_g\omega, \omega\rangle_{L^2\Omega^p(K,g)}$ is non-negative.  This follows directly from Proposition~\ref{prop:qOrbBdry} because the boundary conditions send the boundary terms to zero and $q$ is non-negative.  For similar reasons, $\langle \Delta_g\omega,\eta\rangle_{L^2\Omega^p(K,g)}=\langle \omega, \Delta_g\eta\rangle_{L^2\Omega^p(K,g)}$, i.e., the Laplacian is symmetric on the subspace of $H^2\Omega^p(K,g)$ with either relative or absolute boundary conditions.

We are now in a position to prove the following.

\begin{theorem}
\label{thrm:spectraltheoremNeumannDirichlet}
Let $(K,g)$ be a compact oriented Riemannian orbifold with manifold boundary.
The spectrum of eigenvalues of the $p$-form Laplacian $\Delta_g$
with relative boundary conditions is a discrete set $\{\lambda_{p,k}^D\}_{k=1}^\infty$ with $\lambda_{p,k}^D> 0$
for all $k$. The spectrum of eigenvalues of the $p$-form Laplacian $\Delta_g$ with absolute boundary conditions is
a discrete set $\{\lambda_{p,k}^N\}_{k=0}^\infty$ with $\lambda_{p,k}^N \geq 0$ for all $k$.
In each case, each eigenvalue appears in the spectrum only finitely many times, and the sequences $\{\lambda_{p,k}^D\}$
and $\{\lambda_{p,k}^N\}$ tend to $\infty$ as $k\to \infty$.
\end{theorem}
\begin{proof}
Note that for $\omega\in \HH_D^p(K,g), \Delta_g\omega=0$.  Thus, it suffices to consider the $L^2$-orthogonal complement $\HHDkK^{\obot}$ in $H^1\Omega_D^p(K,g)$.  Using \cite[Corollary 4.2.3]{davies}, if the resolvent operator $(\Delta_g + 1)^{-1}$ is compact, the result will hold. But by \cite[Exercise~4.2 on p.97]{davies}, $(\Delta_g+1)^{-1}$ is compact if and only if $\HHDkK^{\obot}$, with $q$-norm
\begin{equation}
\label{eq:qNorm}
    \Vert \omega\Vert_q := \big(q(\omega,\omega) + \Vert \omega\Vert^2_{L^2\Omega^p(K, g)}\big)^{\frac{1}{2}},
\end{equation}
is compactly embedded in $L^2\Omega^p(K, g)$.  It follows from Proposition~\ref{prop:qH1norms} that the $q$-norm and $H^1$-norm are equivalent on $\HHDkK^{\obot}$.  Thus, by Rellich's theorem, $\HHDkK^{\obot}$ with $H^1$-norm is compactly embedded in $L^2\Omega^p(K,g)$ and the result follows for relative boundary conditions.  By the comment \cite[p.72]{schwarz}, the result holds for absolute boundary conditions as well.
\end{proof}

With this, we can now state the min-max principle as it applies to this particular situation.
By Theorem~\ref{thrm:spectraltheoremNeumannDirichlet} and Theorem~\ref{thrm:Hodge}, the following version of the min-max
principle carries over directly to the orbifold setting.  See \cite[Proposition 7]{JammesPrecrip};
see also \cite[Proposition 3.1]{dodziuk}, \cite[Proposition 2.1]{mcgowan}, and Remark~\ref{rem:ExactCoexactSpec}.

\begin{theorem}[Min-max principle, (co-)exact forms, absolute boundary conditions]
\label{thrm:minmaxNeumannDirichlet}
Let $(K,g)$ be a compact oriented Riemannian orbifold with manifold boundary and let
\[
H^1\Omega^p_{\mathrm{abs}}(K, g) = \big\{\omega\in H^1\Omega^p(K,g) \,\vert\, \bfn\omega= 0, \bfn d\omega=0\big\}.
\]
For $\omega\in H^1\Omega^p_{\mathrm{abs}}(K,g)$, let
\[
    R(\omega)= \frac{\langle d\omega,d\omega\rangle_{L^2\Omega^{p+1}(K,g)}}{\langle\omega,\omega\rangle_{L^2\Omega^p(K,g)}}.
\]
Then for $k\geq 1$,
\[
    \mu_{p,k}^N(K,g)
    = \inf_Y
    \sup_{0\neq d\omega\in Y}R(\omega)
\]
where $Y$ ranges over all $k$-dimensional subspaces of smooth exact $(p+1)$-forms $d\omega\in L^2\Omega^{p}(K,g)$ satisfying absolute boundary conditions.
\end{theorem}

We note too that for the case of functions on $(K,g)$ with Dirichlet conditions, since $\langle \Delta_gf,f\rangle_{L^2(K,g)}=\langle \nabla_gf,\nabla_gf\rangle_{L^2(K,g)}$, the standard min-max formula holds as well.  See \cite[Theorem 4.5.1]{davies}.

\begin{theorem}[Min-max principle, functions, Dirichlet conditions]
\label{thrm:MinMaxFunctionsDirichlet}
Let $(K,g)$ be a compact oriented Riemannian orbifold with manifold boundary $\partial K$. Then
\[
\lambda_{k}^D = \inf_Y\sup_{0\neq f\in Y}\frac{\langle \nabla_gf,\nabla_gf\rangle_{L^2(K,g)}}{\langle f,f\rangle_{L^2(K,g)}}
\]
where $Y$ ranges over all $k$-dimensional subspaces of $H^1_D(K,g)$, the subspace of $H^1(K,g)$ of functions vanishing on $\partial K$.
\end{theorem}

Of course, if $K$ is closed, then Theorems~\ref{thrm:spectraltheoremNeumannDirichlet}, \ref{thrm:minmaxNeumannDirichlet}, and \ref{thrm:MinMaxFunctionsDirichlet} apply to the corresponding $\lambda_{p,k}$ and $\mu_{p,k}$ with vacuous boundary conditions.

We now recall the following lemma, demonstrating the existence of harmonic extensions for orbifolds.

\begin{lemma}
[{\cite[Lemma 3.5]{AMDGHRS}}]
\label{lem:HarmonicExt}
Let $(K,g)$ be a compact oriented Riemannian orbifold with manifold boundary $\partial K$ and let $u \in\mathcal{C}^\infty(\partial K)$.
Then there exists a unique harmonic function $h \in \mathcal{C}^\infty(K)$ such that $h_{\vert \partial K} = u$.
\end{lemma}

The idea of the proof is to use the fact that $K$ can be expressed as the quotient of its orthonormal frame
bundle $\mathcal{F}K$, a smooth manifold, by the almost free action of the orthogonal group $\operatorname{SO}(n)$,
and then $\mathcal{F}K$ admits a metric $\tilde{g}$ with respect to which
\begin{equation}
\label{eq:FrameLaplac}
    \pi^\ast\circ\Delta_{g}
        =   \Delta_{\tilde{g}}\circ\pi^\ast,
        \end{equation}
 where $\pi\co\mathcal{F}K\to K$ is the quotient map.

Using Lemma~\ref{lem:HarmonicExt}, the proof of \cite[Chapter 5, Proposition 1.7]{taylor} with $s = 1/2$ extends directly to the orbifold case, yielding the following.

\begin{lemma}
\label{lem:HarmonicExtSobolev}
Let $(K,g)$ be a compact oriented Riemannian orbifold with nonempty manifold boundary $\partial K$. The map associating to $u \in\mathcal{C}^\infty(\partial K)$ its harmonic extension to $K$ admits a continuous extension $P\co H^{1/2}(\partial K, \partial g) \to H^1(K, g)$. Hence, there is a constant $C$ independent of $u$ such that $\Vert P(u)\Vert_{H^1(K, g)} \leq C\Vert u \Vert_{H^{1/2}(\partial K, \partial g)}$.
\end{lemma}

Using the same method as Lemma~\ref{lem:HarmonicExt}, we can generalize
the Hopf Maximum Principle \cite[Chapter 5, Propositon 2.3]{taylor} to the case of orbifolds.

\begin{theorem}[The Hopf Maximum Principle for Riemannian orbifolds]
\label{thrm:HopfMaxPrin}
Suppose $U$ is a connected open subset with nonempty manifold boundary $\partial U$
of a compact oriented Riemannian orbifold $(K, g)$. If $f\in\mathcal{C}^2(U)$
is continuous on the closure $\overline{U}$ and $\Delta_g f \geq 0$, then either
$f$ is constant or
\begin{equation}
\label{eq:HopfMaxPrin}
    f(x) < \sup\limits_{y\in\partial U} f(y)
    \quad\quad\forall x\in U.
\end{equation}
\end{theorem}
\begin{proof}
Given $f$ as in the statement, it is clear that $\pi^\ast f \in\mathcal{C}^2(\mathcal{F}U)$
is continuous on the closure of $\mathcal{F}U$. Then by Equation~\eqref{eq:FrameLaplac},
$\Delta_g f \geq 0$ implies that $\Delta_{\tilde{g}}\circ\pi^\ast f \geq 0$. Applying
the maximal principal \cite[Chapter 5, Propositon 2.3]{taylor} for manifolds to $\pi^\ast f$ completes the
proof.
\end{proof}


\subsection{The Laplacian on orbifolds with connected sum metrics}
\label{subsec:LaplacianConnectedSum}

Suppose that $O$ is a closed oriented orbifold with connected sum metric $g$ as in Definition~\ref{def:ConnSum}.  We define a bilinear form $q$ on $H^1\Omega^p(O,g)$ using Equation~\eqref{eqn:qform}.  Namely, for $\omega=(\omega_1,\omega_2)$ and $\eta=(\eta_1,\eta_2)$ in $H^1\Omega^p(O,g)$, let
\[
q(\omega,\eta)=q_1(\omega_1,\eta_1) + q_2(\omega_2,\eta_2)
\]
where $q_i$ denotes the bilinear form on $(O_i,g_i)$, $i=1,2$ respectively.

Similarly, for $\omega=(\omega_1,\omega_2)\in H^2\Omega^p(O,g)$, define the Laplacian by
\[
\Delta_g\omega = (\Delta_{g_1}\omega_1, \Delta_{g_2}\omega_2).
\]

We first note that as with orbifolds with boundary, for orbifolds with connected sum metrics, the Laplacian is induced by $q$.

\begin{proposition}
\label{prop:QveLaplacian}
Let $O$ be a closed oriented orbifold with connected sum metric $g$.
For $\omega=(\omega_1,\omega_2)\in H^2\Omega^p(O, g)$ and $\eta=(\eta_1,\eta_2)\in H^1\Omega^p(O, g)$,
the bilinear form $q$ induces the Laplacian, i.e.,
\[
    q(\omega,\eta) =\langle \Delta_g \omega,\eta \rangle_{L^2\Omega^p(O, g)}.
\]
\end{proposition}
\begin{proof}
Because $\omega\in H^2\Omega^p(O,g)$, we have $\bft_1\omega_1=\bft_2\omega_2$, $\star\bfn_1\omega_1=-\star\bfn_2\omega_2$, $\bft\delta\omega_1=\bft_2\delta\omega_2$, and $\star\bfn_1d\omega_2=-\star\bfn_2d\omega_2$.  Furthermore, since $\eta\in H^1\Omega^p(O,g)$, $\bft_1\eta_1=\bft_2\eta_2$ and $\star\bfn_1\eta_1=-\star\bfn_2\eta_2$.

The result then follows directly from the definition of $q$ and Proposition~\ref{prop:qOrbBdry}.
\end{proof}

Furthermore, using a similar analysis and the fact that $q$ is a symmetric, non-negative bilinear form, we see directly that for $\omega, \eta\in H^2\Omega^p(O,g)$, $\langle\Delta_g\omega,\omega\rangle_{L^2\Omega^p(O,g)}$ is non-negative and symmetric on $H^2\Omega^p(O,g)$.

We thus have the following, in analogy with Theorem~\ref{thrm:spectraltheoremNeumannDirichlet}.

\begin{theorem}
\label{thrm:spectraltheorem}
Let $O$ be a closed oriented $n$-dimensional orbifold with connected sum metric $g$.
For $p=0,\ldots,n$, the spectrum of eigenvalues of the Laplacian $\Delta_g$ acting on $p$-forms
is a discrete set $\{\lambda_{p,k}\}$ with $\lambda_{p,k}\geq 0$ for all $k$.
Each eigenvalue appears in the spectrum only finitely many times, and the sequence $\{\lambda_{p,k}\}$ tends to $\infty$ as
$k\to \infty$. There is also a complete orthonormal basis of $L^2\Omega^p(O, g)$ consisting of eigenforms
$\{\omega_{p,k}\}_{k=0}^{\infty}$ of $\Delta_g$.
\end{theorem}

\begin{proof}
The claim follows from the equivalence of the $q$-norm of Equation~\eqref{eq:qNorm} and the norm $H^1$-norm. To show that these two norms are equivalent, one uses the proof of Theorem~\ref{thrm:spectraltheoremNeumannDirichlet} in the case that the form has support away from the boundary. In the case that the support is included in a neighborhood of the boundary, one uses instead
\cite[Proposition~1.2]{annecolbois}, which readily generalizes to the orbifold case with manifold boundary as the argument is localized near the boundary.
\end{proof}


\section{Approximating the spectrum on functions}
\label{sec:SpecDefDim3Func}

In this section, we consider the convergence of eigenvalues of the Laplacian acting on functions in two ways.  The first is a generalization of work of Colin de Verdi\`ere. We begin by reviewing his results which detail when two quadratic forms have close finite spectrum. We then generalize to the orbifold setting his result that allows us to perturb a given orbifold metric in order to approximate the finite spectrum of the Laplacian with that of the Neumann problem on a domain in the orbifold. The second is a generalization of work originally due to Rauch and Taylor.  In this method, we show that we may remove a ball from an orbifold in such a way that as the radius of the ball shrinks to zero, the spectrum of the Neumann problem on the resulting orbifold domain tends to the Laplace spectrum of the orbifold, with its original metric.  This generalizes results known to hold for manifolds; see relevant references below.


\subsection{The defect of quadratic forms}
\label{subsec:PrelimQuadForms}

Recall the following.

\begin{definition}[{\cite[Definition I.2]{VerdierMultiplicite}}]
\label{def:N-spectral-defect}
Let  $(E_0, E_0, q_0)$ and $(E_1, E_1, q_1)$ be two positive quadratic forms defined on the finite-dimensional
Euclidean spaces $E_0$ and $E_1$ with inner products $\langle\cdot,\cdot\rangle_0$ and
$\langle\cdot,\cdot\rangle_1$ respectively. We say that $q_0$ and $q_1$ have \emph{defect $\leq \ve$}
if there is an isometry $U_{0,1} \co E_0 \to E_1$ such that
\[
    \| q_1\circ  U_{0, 1} - q_0 \| \leq \ve,
\]
where the norm $\|\cdot\|$ is that defined in Equation~\eqref{eq:QuadFormNorm}.

Let $(\mathcal{H}_i, D(q_i), q_i)$ be quadratic forms for $i=0,1$ and fix
$N \in \N$. For each $i$, let $E_i$ be the sum of the eigenspaces associated to the first $N$ nonzero eigenvalues of $q_i$. We say that
$(\mathcal{H}_0, D(q_0), q_0)$ and $(\mathcal{H}_1, D(q_1), q_1)$ have
\emph{$N$-spectral defect $\leq \ve$} if $\tilde{q}_0 := q_{0 \vert E_0}$ and $\tilde{q}_1 := q_{1 \vert E_1}$ have defect $\leq\ve$.  A small $N$-spectral defect means that the first $N$ eigenvalues of $q_0$ are close to the first $N$ eigenvalues of $q_1$.
\end{definition}

The following simple observation concerning approximate transitivity of the $N$-spectral defect will be useful in the sequel.

\begin{lemma}
\label{lem:NSpecDefTransitive}
Suppose $N \in \N$, $(\mathcal{H}_i, D(q_i), q_i)$ are positive quadratic forms for $i=0,1,2$,
$E_i$ is the sum of the eigenspaces associated to the first $N$ eigenvalues of $q_i$, and
$\tilde{q}_i := q_{i \vert E_i}$. If $(\mathcal{H}_0, D(q_0), q_0)$ and $(\mathcal{H}_1, D(q_1), q_1)$ have
$N$-spectral defect $\leq\ve$ and $(\mathcal{H}_1, D(q_1), q_1)$ and $(\mathcal{H}_2, D(q_2), q_2)$
have $N$-spectral defect $\leq\ve$, then $(\mathcal{H}_0, D(q_0), q_0)$ and $(\mathcal{H}_2, D(q_2), q_2)$
have $N$-spectral defect $\leq2\ve$.
\end{lemma}
\begin{proof}
Let $U_{0,2}:= U_{1,2}\circ U_{0,1}$. Then
\[
    \| \tilde{q}_2\circ U_{0, 2} - \tilde{q}_0 \|
    \leq   \| \tilde{q}_2\circ U_{1,2}\circ U_{0,1} - \tilde{q}_1\circ U_{0,1}\| + \| \tilde{q}_1\circ U_{0,1} - \tilde{q}_0 \| \leq 2\ve.
    \qedhere
\]
\end{proof}

We say that a quadratic form $(\mathcal{H}, D(q), q)$ satisfies hypothesis \eqref{eq:HypAst} if for some fixed $M, N$, and $\delta>0$, the nonzero eigenvalues $\lambda_i$ of $q$ satisfy
\begin{equation}
\label{eq:HypAst}
    \tag{$\ast$}
    0 < \lambda_1 \leq \lambda_2 \leq \cdots \leq \lambda_N < \lambda_N + \delta \leq \lambda_{N+1} \leq M.
\end{equation}
See \cite[page 257]{VerdierMultiplicite}.

Recall the following.

\begin{theorem}[{\cite[Theorem I.7]{VerdierMultiplicite} and \cite[Lemme 16]{JammesPrecrip}}]
\label{thrm:1.7-Verdier-first}
Let $\eta>0$ and $q$ be a positive quadratic form  on $\mathcal{H}$ with domain
$\mathcal{D}(q) = \mathcal{K}_0 \oplus \mathcal{K}_\infty$ where $\mathcal{K}_0$ and $\mathcal{K}_\infty$
are $q$-orthogonal subspaces. Suppose that for some $\delta$, $M$, and $N$, hypothesis \eqref{eq:HypAst} is satisfied for
$q_0 := q_{\vert\mathcal{K}_0 }$ and that, for $x \in \mathcal{K}_\infty$, we have:
\[
    q(x) \geq C \Vert x \Vert^2
\]
for some constant $C$ depending only on $\delta,M,N$, and $\eta$. Then the forms $q_0$ and $q$ have $N$-spectral
defect $\leq \eta$.
\end{theorem}

\begin{theorem}
[{\cite[Theorem 1.8 and page 266]{VerdierMultiplicite} and \cite[Lemme 17 and  Remarque 18]{JammesPrecrip}}]
\label{thrm:Thm-1.8-Verdier}
Let $\mathcal{H}$ be a Hilbert space with norm $\Vert \cdot \Vert$. Let $\eta>0$, and
let $q$ and $q_m$ for each $m \in \N$ be positive quadratic forms on $\mathcal{H}$
with common domain $D(q)$. Let $\Vert\cdot\Vert_m$ be additional norms on $\mathcal{H}$,
and assume that:
\begin{enumerate}
\item   there are constants $C_1, C_2 >0$ such that
\[
    C_1\, \Vert x \Vert \leq \Vert x\Vert_m  \leq C_2\,  \Vert x \Vert + \ve_m q(x)^{1/2}
    \; (\text{with } \ve_m\to 0)
    \quad\forall x \in \mathcal{H}.
\]
\item for all $x \in D(q)$, we have
\[
    \lim_{m \to \infty} \Vert x\Vert_m=\Vert x\Vert.
\]
\item for all $x \in D(q)$, we have
\[
    q(x) \leq  q_m(x), \quad \forall m \in \N.
\]
\item for all $x \in D(q)$, we have
\[
    \lim_{m \to \infty} q_m(x) = q(x).
\]
\end{enumerate}
If in addition $q$ satisfies hypothesis \eqref{eq:HypAst} for some $\delta$, $M$, and $N$, then there exists an
$m_0$ such that for $m \geq m_0$, $q_m$ and $q$ have $N$-spectral defect $\leq \eta$.
\end{theorem}


\subsection{The approach of Colin de Verdi\`ere}
\label{subsec:SpecDefVerdier}

Our first main result is the generalization of \cite[Theorem III.1]{VerdierMultiplicite} to the case
of Riemannian orbifolds.  Our proof follows that of \cite[Theorem 11]{JammesPrecrip}.  We note that we make a minor change to the definition of $g_{\ve}$ from \cite{VerdierMultiplicite}, replacing $\ve$ with $\ve^2$, in order to be consistent with \cite{JammesPrecrip} and other arguments in this paper. Note that by Theorem~\ref{thrm:spectraltheoremNeumannDirichlet}, for any $N^\prime > 0$, there is an $N > N^\prime$, a $\delta > 0$, and an $M > 0$ such that hypothesis \eqref{eq:HypAst} is satisfied for the eigenvalues of the Neumann problem on $\overline{U}$.

\begin{theorem}
\label{thrm:SpecDefDim3}
Let $(O, g)$ be a closed oriented Riemannian orbifold of dimension $n \geq 3$ and let $\overline{U}\subset O$
be a compact orbifold domain with interior $U$ and manifold boundary
$\Theta:=\partial\overline{U}=\partial({O\smallsetminus U})$. Assume that for some $\delta$, $M$, and $N$,
hypothesis \eqref{eq:HypAst}
holds for the eigenvalues of the Neumann problem on $\overline{U}$.
Then for any  $\eta > 0$  there exists a $\mathcal{C}^\infty$ metric $h_\eta$ on $O$ such that
$h_{\eta\vert\overline{U}} = g_{\vert\overline{U}}$ and such that the $N$-spectral defect of the
Neumann problem on $\overline{U}$ and the Laplacian on  $(O, h_\eta)$  is $\leq \eta$.
\end{theorem}
\begin{proof}
We first give a brief outline of the proof.
To begin, we will define a piecewise smooth metric $g_\ve$ given by shrinking $g$ on $O\smallsetminus \overline{U}$
and approximate $g_\ve$ with smooth metrics $g_{m,\ve}$. The first step will then be an application of
Theorem~\ref{thrm:Thm-1.8-Verdier} to show that for quadratic forms $q_\ve$ and $q_{m,\ve}$ on $H^1(O,g)$
associated to $g_\ve$ and $g_{m,\ve}$, respectively, the $N$-spectral defect can be made arbitrarily small for
$N$ large. In the second step, we give a decomposition $\mathcal{K}_0\oplus\mathcal{K}_\infty$
of $H^1(O,g)$ and apply Theorem~\ref{thrm:1.7-Verdier-first} to show that for small $\ve$, the
$N$-spectral defect of $q_\ve$ and $q_{\ve\vert\mathcal{K}_0}$ can be made arbitrarily small.
The third step is to show that $\ve$ can be chosen so that $q_{\ve\vert\mathcal{K}_0}$ and the standard quadratic form associated to
the Neumann problem on $(\overline{U},g)$ have arbitrarily small $N$-spectral defect, also
an application of Theorem~\ref{thrm:1.7-Verdier-first}. Then an application of Lemma~\ref{lem:NSpecDefTransitive}
demonstrates that, for sufficiently large $N$ and small $\ve$, the $N$-spectral defect of the Neumann problem on
$(\overline{U},g)$ and $(O, g_{m,\ve})$ can be made arbitrarily small.

Fix $\ve >0$ and choose a decreasing sequence $F_m^\ve \geq F_{m+1}^\ve$ of positive smooth functions on $O$ such that $F^\ve_{m\vert\overline{U}} = 1$ for each $m$, $\ve\leq F_m^\ve(x) \leq 1$ for each $m$ and $x \in O \smallsetminus \overline{U}$, and $\lim_{m\to\infty} F_m^\ve(x) = \ve$ for each $x \in O \smallsetminus \overline{U}$. Define the metrics $g_\ve, g_{m,\ve}$ on $O$ by
\begin{equation}
\label{eq:Defgve}
    g_\ve = \begin{cases}
    g       &\mbox{on } \;      \overline{U},    \\
    \ve^2 g &\mbox{on } \;      O \smallsetminus U,
\end{cases}\quad
    g_{m,\ve} = (F_m^\ve)^2 g.
\end{equation}
Note that $g_\ve$ is not a connected sum metric in the sense of Definition~\ref{def:ConnSum}, but rather
is multi-valued on the set $\Theta$, which has measure zero, and smooth on $\overline{U}$ and $O \smallsetminus U$
separately. Below, integrals on
$\overline{U}$ and $O \smallsetminus U$ will be with respect to the definition of $g_\ve$ that is continuous on that piece.
Moreover, the smooth metrics $g_{m, \ve} := (F_m^\ve)^2 g$ converge to the discontinuous metric $g_\ve$
on $O\smallsetminus\Theta$.

Let $\mathcal{H}=L^2(O,g)$ and let $D(q)=H^1(O,g)$, which will be the common domain of the quadratic forms we now
define. We will use $(\varphi_1, \varphi_2)$ to denote an element
$\varphi\in D(q)$ where $\varphi_1 = \varphi_{|\overline{U}}$ and $\varphi_2 = \varphi_{|O\smallsetminus U}$.
For metrics $g_{\ve}$ and $g_{m,\ve}$ on $O$, define quadratic forms $q_{\ve}$ and $q_{m,\ve}$ on $D(q)$ by
\begin{align*}
q_{\ve}(\varphi)&=\Vert \nabla_{g_{\ve}}\varphi\Vert_{L^2(O,g_{\ve})}^2\\
&=\int_{\overline{U}} \vert \nabla_{g_\ve}\varphi_1\vert_{g_\ve}^2\, dv_{g_\ve} +\int_{O\smallsetminus U}\vert \nabla_{g_{\ve}}\varphi_2\vert_{g_\ve}^2\, dv_{g_{\ve}}\\
&=\int_{\overline{U}} \vert \nabla_g\varphi_1\vert_g^2\, dv_g + \ve^{n-2}\int_{O\smallsetminus U}\vert \nabla_g\varphi_2\vert_g^2\,dv_g,
\end{align*}
where the notation $\vert \cdot\vert_g$ indicates that the norm is taken using the metric $g$, and
\begin{align*}
q_{m,\ve}(\varphi)=\int_{\overline{U}} \vert \nabla_g\varphi_1\vert_g^2\, dv_g +\int_{O\smallsetminus U}(F_m^{\ve})^{n-2}\vert \nabla_g\varphi_2\vert_g^2\, dv_{g}.
\end{align*}

Define the norms $\Vert\cdot\Vert_{\ve}$ and $\Vert\cdot\Vert_{m,\ve}$ by
\begin{align*}
\Vert\varphi\Vert_{\ve}^2 &=\Vert\varphi\Vert_{L^2(O,g_{\ve})}^2\\
&=\int_{\overline{U}} \vert \varphi_1\vert^2\, dv_{g_{\ve}} +\int_{O\smallsetminus U}\vert \varphi_2\vert^2\, dv_{g_{\ve}}\\
&=\int_{\overline{U} }\vert \varphi_1\vert^2\, dv_g + \ve^{n}\int_{O\smallsetminus U}\vert \varphi_2\vert^2\,dv_g
\end{align*}
and
\begin{align*}
\Vert \varphi\Vert_{m,\ve}^2&=\Vert\varphi\Vert_{L^2(O,g_{m,\ve})}^2\\
&=\int_{\overline{U} }\vert \varphi_1\vert^2\, dv_g +\int_{O\smallsetminus U}(F_m^{\ve})^{n}\vert \varphi_2\vert^2\, dv_{g}.
\end{align*}

Since $\ve\leq F_m^\ve(x) \leq 1$ on $O\smallsetminus \overline{U}$, we have that $q_\ve(\varphi) \leq q_{m,\ve}(\varphi)$ for all $\varphi\in D(q)$.  Recalling that $\lim_{m\to\infty} F_m^\ve(x) = \ve$ on $O \smallsetminus \overline{U}$, it follows that $g_{m,\ve}\to g_{\ve}$ pointwise on $O\smallsetminus \overline{U}$, and hence $\lim_{m\to\infty} q_{m,\ve}(\varphi) = q_{\ve}(\varphi)$ and $\lim_{m\to\infty}\Vert\varphi\Vert_{m,\ve}=\Vert\varphi\Vert_{\ve}$ for all $\varphi\in D(q)$.  Finally,
\[
    \int_{O\smallsetminus U}|\varphi_2|^2 \, dv_{g_{\ve}}
    \leq
    \int_{O\smallsetminus U}|\varphi_2|^2 \, dv_{g_{m,\ve}}
    \leq
    \int_{O\smallsetminus U}|\varphi_2|^2 \, dv_{g}
    =
 \frac{1}{\ve^{n}}
    \int_{O\smallsetminus U}|\varphi_2|^2 \, dv_{g_{\ve}},
\]
implying that $\Vert\varphi\Vert_\ve^2 \leq \Vert\varphi\Vert_{m,\ve}^2 \leq \ve^{-n} \Vert\varphi\Vert_\ve^2$.

It follows that the $q_{\ve}$ and $q_{m,\ve}$ satisfy the hypotheses of
Theorem~\ref{thrm:Thm-1.8-Verdier} and hence that the $N$-spectral defect of the quadratic form
$q_\ve$ on $D(q)=H^1(O,g)$ associated to the metric $g_{\ve}$ and the quadratic form $q_{m,\ve}$
associated to the metric $g_{m,\ve}$ can be made arbitrarily small for $m$ large. This completes the first step.

We next will apply Theorem~\ref{thrm:1.7-Verdier-first} for suitable subspaces
$\mathcal{K}_0$ and $\mathcal{K}_\infty$ of $D(q)=H^1(O,g)$ in order to show that for small $\ve$, $q_{\ve}$ and $q_{\ve,0}=q_{\ve\vert\mathcal{K}_0}$ have $N$-spectral defect arbitrarily small.
To define these spaces, first note that for $\varphi\in H^1(O,g)$, the restriction $\varphi_{\vert\Theta} \in H^{1/2}(\Theta,\partial g)$ by the Trace theorem, Theorem~\ref{thrm:trace}~(2). Let $P\co H^{1/2}(\Theta, \partial g) \to H^1(O\smallsetminus U, g)$ denote the harmonic extension operator; see Lemmas~\ref{lem:HarmonicExt} and \ref{lem:HarmonicExtSobolev}, and note that the harmonic extensions with respect to $g$ and $g_\ve$ coincide.
For $\varphi=(\varphi_1,\varphi_2)\in D(q)=H^1(O,g)$, we can express
\begin{equation}
\label{eq:HarmExt}
    (\varphi_1, \varphi_2) = \big(\varphi_1, P(\varphi_{1\vert\Theta})\big)
                    + \big(0, \varphi_2- P(\varphi_{1\vert\Theta})\big).
\end{equation}
Note that since $(\varphi_1,\varphi_2)\in H^1(O,g)$, so are both $(\varphi_1,P(\varphi_{1\vert\Theta}))$ and $(0,\varphi_2-P(\varphi_{1\vert\Theta}))$.

With this, we define
\[
    \mathcal{K}_0
        :=  \big\{ (\varphi_1, P(\varphi_{1\vert\Theta}))\big\}\subset D(q)
\]
and
\[
    \mathcal{K}_\infty
        :=  \big\{ (0, \psi_2)\big\} \subset D(q).
\]
Note that for $(0,\psi_2)\in \mathcal{K}_{\infty}$, $\psi_{2\vert\Theta}=0|_{\Theta}$ in $L^2(\Theta,\partial g)$.   By Equation~\eqref{eq:HarmExt}, it is clear that $D(q)=\mathcal{K}_0\oplus\mathcal{K}_{\infty}$.  Furthermore, $\mathcal{K}_0$ and $\mathcal{K}_{\infty}$ are $q_{\ve}$-orthogonal.   Indeed, for $\varphi,\psi\in D(q)=H^1(O,g)$, the symmetric bilinear form $\hat q_{\ve}$ associated to $q_{\ve}$ is given by
\[
    \hat q_{\ve}(\varphi,\psi)=\frac{1}{2}\big( q_{\ve}(\varphi +\psi) -q_{\ve}(\varphi)-q_{\ve}(\psi)\big).
\]
By direct computation, we find using Green's Theorem (Theorem~\ref{thrm:Green}) that
\begin{align*}
    2\hat q_{\ve} \Big(\big(\varphi_1,P(\varphi_{1\vert\Theta}),(0,\psi_2)\big)\Big)
        &=  \int_{O\smallsetminus U}\langle\Delta_{g_{\ve}}P(\varphi_{1\vert\Theta}),\psi_2\rangle_{g_\ve}\, dv_{g_{\ve}}
            - \int_{\Theta} \psi_2\Big(\nu\big( P(\varphi_{1\vert\Theta})\big)\Big) dv_{\partial g_\ve} \\
        &\quad  +\int_{O\smallsetminus U}\langle\Delta_{g_{\ve}}\psi_2, P(\varphi_{1\vert\Theta})\rangle_{g_\ve}\, dv_{g_{\ve}}
                -\int_{\Theta} P(\varphi_{1\vert\Theta})\big(\nu( \psi_2 )\big) dv_{\partial g_\ve}.
\end{align*}
The first term vanishes because $P(\varphi_{1\vert\Theta})$ is harmonic; the second because $\psi_2$ is zero on $\Theta$;
the third term vanishes by the Hodge decomposition theorem, Theorem~\ref{thrm:Hodge}, as
$\Delta_{g_{\ve}}\psi_2\in\mathcal{E}^0(O\smallsetminus U,g_\ve)$ and hence is orthogonal to
$P(\varphi_{1\vert\Theta})$; and the fourth term vanishes because $\nu( \psi_2 ) = \nu(0) = 0$ by the definition
of $H^1(O,g)$. Therefore, $\hat q_{\ve} = 0$.

Now, let $C= {\ve^{-2} \lambda^D_1(O\smallsetminus U)},$  where $\lambda^D_1(O\smallsetminus U)$
is the first nonzero eigenvalue of the Dirichlet problem on $O\smallsetminus U$ with respect to the
metric $g$. Recall that $\nabla_g=(d + \delta)$ (with $\delta = 0$ on functions) and $\nabla_g^2=\Delta_g$. Then by Theorem~\ref{thrm:MinMaxFunctionsDirichlet}, $\lambda_1^D(O\smallsetminus U)$ is the infimum of the Rayleigh quotient
\[
    \frac{\Vert \nabla_g \psi_2 \Vert_{L^2(O\smallsetminus U, g)}}{\Vert \psi_2 \Vert_{L^2(O\smallsetminus U, g)}}
\]
where $\psi_2$ ranges over nonzero $\mathcal{C}^1$ functions on $O\smallsetminus U$ that
vanish on the boundary $\Theta$.  Thus, we have that for $\psi=(0,\psi_2)\in\mathcal{K}_\infty$,
\begin{align*}
    q_\ve(\psi)
        &=          {\ve}^{n-2}  \int_{O\smallsetminus U} \vert\nabla_g \psi_2 \vert_g^2 \, dv_g
        \\&\geq     {\ve}^{n-2} \lambda^D_1(O\smallsetminus U) \int_{O\smallsetminus U} \vert\psi_2 \vert^2 \, dv_g
        \\&=        \ve^{-2} \lambda^D_1 (O\smallsetminus U)
                    \Big(\ve^n \int_{O\smallsetminus U} \vert\psi_2 \vert^2 \, dv_g\Big)
        =           C  \Vert \psi \Vert_\ve^2.
\end{align*}
As hypothesis \eqref{eq:HypAst} is satisfied for the eigenvalues
of the Neumann problem on $\overline{U}$, it follows by Theorem~\ref{thrm:1.7-Verdier-first} that
the quadratic forms $q_\ve$ on $D(q)$ and $q_{\ve\vert\mathcal{K}_0}$ on $\mathcal{K}_0$
(both with respect to the norm $\Vert\cdot\Vert_\ve$)
have $N$-spectral defect arbitrarily small for small enough $\ve$.
This completes the second step.

We now turn our attention to the third step.  First, note that $\mathcal{K}_0$ can be identified with $H^1(\overline{U},g)$ via the map $(\varphi_1,P(\varphi_{1\vert\Theta}))\mapsto\varphi_1$. Define the quadratic form
\[
    q(\varphi)  := \int_{\overline{U}}  |\nabla_g\varphi|_g^2 \, dv_g
\]
on $H^1(\overline{U},g)$ and the norm
\[
    \Vert\varphi \Vert_0^2 :=\Vert\varphi\Vert_{L^2(\overline{U}, g)}^2=  \int_{\overline{U}} |\varphi|^2 \, dv_g.
\]
Then $q$ is the quadratic form associated to the Neumann problem on $\overline{U}$ with domain
$H^1(\overline{U},g)$, see Proposition~\ref{prop:qOrbBdry} and \cite[Section 3.1]{Anne-Post}. Suppose that $\varphi=(\varphi_1,P(\varphi_{1\vert\Theta}))\in\mathcal{K}_0$.  From Equation~\eqref{eq:HarmExt},
it follows that the restriction $q_{\ve,0}$ of $q_{\ve}$ to ${\mathcal{K}_0}$ is given by
\[
    q_{\ve,0}(\varphi)  = \int_{\overline{U}} \vert\nabla_g\varphi_1\vert_g^2 \, dv_g
        +  \ve^{n-2} \int_{O\smallsetminus U} \vert\nabla_g P(\varphi_{1\vert\Theta}) \vert_g^2 \, dv_{g},
\]
and similarly the restriction $\Vert \cdot\Vert_{\ve,0}$  of $\Vert \cdot\Vert_{\ve}$ to $\mathcal{K}_0$
is given by:
\[
    \Vert\varphi \Vert_{{\ve},0}^2 := \Vert \varphi\Vert_{L^2(O,g_{\ve})}^2
        = \int_{\overline{U}} \vert\varphi_1\vert^2 \, dv_{g} + \ve^n\int_{O\smallsetminus U} \vert P(\varphi_{1\vert\Theta}) \vert^2 \, dv_{g}.
\]

To apply Theorem~\ref{thrm:Thm-1.8-Verdier}, fix a sequence $\epsilon_m \to 0$ and consider $q_{m,0}:= q_{\epsilon_m,0}$ and norms $\Vert\cdot\Vert_{m,0}:=\Vert\cdot\Vert_{\epsilon_m,0}$.  We note that hypotheses (2), (3), and (4) of the theorem are obviously satisfied, as is the first inequality of (1). To verify the second inequality of (1), let $\varphi=(\varphi_1,P(\varphi_{1\vert\Theta}))\in\mathcal{K}_0$. By Lemma~\ref{lem:HarmonicExtSobolev}, there is a constant $C > 0$ independent of $\varphi$ such that
$\Vert P(\varphi_{1\vert\Theta})\Vert_{H^1(O\smallsetminus U, g)} \leq C\Vert \varphi_{1\vert\Theta} \Vert_{H^{1/2}(\Theta, \partial g)}$, i.e.,
\[
    \int_{O\smallsetminus U} \vert P(\varphi_{1\vert\Theta})\vert^2 \, dv_g \leq
    \Vert P(\varphi_{1\vert\Theta})\Vert_{H^1(O\smallsetminus U, g)}^2
        \leq    C^2 \Vert \varphi_{1\vert\Theta} \Vert_{H^{1/2}(\Theta, \partial g)}^2.
\]
Then by Theorem~\ref{thrm:trace}~(2), there is a constant $C^\prime > 0$ independent of $\varphi$ such that
\[
    \Vert \varphi_{1\vert\Theta} \Vert_{H^{1/2}(\Theta, \partial g)}
        \leq C^\prime \Vert \varphi_1 \Vert_{H^1(\overline{U}, g)}.
\]
Finally, by Gaffney's inequality, Lemma~\ref{lem:Gaffney}, as
$\Vert d\varphi\Vert_{L^2\Omega^1(\overline{U},g)}^2 = q(\varphi)$ and $\delta\varphi = 0$, there is a
$C^{\prime\prime} > 0$ independent of $\varphi$ such that
\[
    \Vert \varphi_1 \Vert_{H^1(\overline{U}, g)}
        \leq C^{\prime\prime} \big( \Vert\varphi_1\Vert_{L^2(\overline{U}, g)} + q(\varphi)^{1/2} \big).
\]
Combining these observations and setting $C^{\prime\prime\prime} = C C^\prime C^{\prime\prime}$,
we have
\[
    \int_{O\smallsetminus U} \vert P(\varphi_{1\vert\Theta})\vert^2 \, dv_g
        \leq (C^{\prime\prime\prime})^2 \big( \Vert\varphi_1\Vert_{L^2(\overline{U}, g)} + q(\varphi)^{1/2} \big)^2,
\]
and hence, by a straightforward computation,
\[
    \Vert\varphi \Vert_0
        \leq    \Vert\varphi \Vert_{m,0}
        \leq    \Vert\varphi_1\Vert_{L^2(\overline{U}, g)} \big( 1+  \epsilon^{n/2}_m  C^{\prime\prime\prime} \big)
                + \epsilon^{n/2}_m  C^{\prime\prime\prime}   q(\varphi)^{1/2} .
\]
This demonstrates that hypothesis (1) of Theorem~\ref{thrm:Thm-1.8-Verdier} is satisfied, and we conclude
that the Neumann quadratic form $q$ on $\overline{U}$ and the form $q_{\ve}$ restricted to
${\mathcal{K}_0}$ have arbitrarily small $N$-spectral defect for sufficiently small $\ve$. This completes
the third step.

With this, we conclude that, for sufficiently large $m$ and small $\ve$, the $N$-spectral defect of $q_\ve$
and the Laplacian associated to $g_{m,\ve}$ can be made $\leq \eta/3$, the $N$-spectral defect of $q_\ve$
and $q_{\ve\vert\mathcal{K}_0^\ve}$ is $\leq \eta/3$, and the $N$-spectral defect of $q_{\ve\vert\mathcal{K}_0^\ve}$
and the standard quadratic form associated to the Neumann problem on $\overline{U}$ $\leq\eta/3$.
By Lemma~\ref{lem:NSpecDefTransitive}, letting $h_\eta:= g_{m,\ve}$, it follows that the $N$-spectral
defect of the Laplacian of $(O, h_\eta)$ and the Neumann problem on $\overline{U}$ is $\leq\eta$, completing
the proof.
\end{proof}

\begin{remark}
\label{rem:FuncDim2}
For manifolds, Theorem~\ref{thrm:SpecDefDim3} has recently been extended to the case of dimension $n = 2$ in
\cite[Theorem 2.1]{colboiselsufi} using similar techniques to those above.
\end{remark}


\subsection{The approach of Rauch and Taylor}
\label{subsec:SpecDefRauchTaylor}

We now consider an alternate method to approach the question of convergence of eigenvalues following the work of
Rauch and Taylor \cite{rauchtaylor}; see also \cite{anneperturbNeumann,annespectre,annecolboisoperateur,Anne-Post,chavel-feldman-domains,chavel-feldman-less}.

In~\cite{rauchtaylor}, Rauch and Taylor proved that when we remove a ball of radius $\ve$ from a Euclidean domain and endow the newly-formed boundary with Dirichlet or Neumann conditions, the eigenvalues of the resulting $\ve$-Laplacian acting on functions tend, as $\ve\to 0$, to the eigenvalues of the Laplacian on the original domain. This result has been extended to removing small open balls in compact, connected manifolds \cite{chavel-feldman-domains,chavel-feldman-less}, tubular neighborhoods of submanifolds \cite{anneperturbNeumann}, the Laplacian acting on differential forms in these settings \cite{annespectre,annecolboisoperateur}, and non-compact complete manifolds with infinitely many balls removed \cite{Anne-Post}. Here, we are interested in the following.

Let $(O,g)$ be a closed oriented Riemannian orbifold with Laplace eigenvalues $\lambda_k(O,g)$. Suppose that $p\in O$ is an isolated singular point and $\ve > 0$ is small enough so that $\overline{B(p,\ve)}$ contains no other singularities. Let $U_\ve = O\smallsetminus B(p,\ve)$, and then we obtain eigenvalues $\lambda_k^N(U_\ve, g_{\vert U_{\ve}})$ of the Laplacian acting on functions on $U_\ve$ with Neumann boundary conditions.

\begin{theorem}
\label{thrm:RauchTaylor}
Let $(O,g)$ be a closed oriented Riemannian orbifold of dimension $n \geq 2$. With notation as above, we have for $k=0,1,2,\dots$,
\[
    \lim_{\ve\to 0}\lambda_k^N(U_\ve, g_{\vert U_{\ve}})=\lambda_k(O,g).
\]
\end{theorem}
\begin{proof}
This result is proven by a direct generalization of the corresponding result for manifolds given in \cite[Section 2]{anneperturbNeumann}. First, one establishes the corresponding result for Dirichlet eigenvalues,
$\lim_{\ve\to 0}\lambda_k^D(U_\ve, g_{\vert U_{\ve}})=\lambda_k(O,g)$, see also \cite[Proposition (N)]{chavel-feldman-domains}, implying that $\limsup_{\ve\to 0}\lambda_k^N(U_\ve, g_{\vert U_{\ve}}) \leq \lambda_k(O,g)$. For the opposite inequality, Ann\'{e} defines a family of extension operators $P_\ve\co H^1(U_\ve, g_{\vert U_{\ve}})\to H^1(O,g)$ given by harmonic extension from $\partial B(p,\ve)$ to $\overline{B(p,\ve)}$, see Lemma~\ref{lem:HarmonicExtSobolev} and \cite[Definition 4.1 (iii)]{Anne-Post}, and demonstrates that $\| P_\ve\|$ is bounded independently of $\ve$. It follows that for a sequence $\ve_p\to 0$, if $\varphi_k(\ve_p)\in H^1(U_\ve, g_{\vert U_{\ve}})$ are orthonormal eigenfunctions corresponding to the $\lambda_k^N(U_\ve, g_{\vert U_{\ve}})$, there is a subsequence such that
$P_{\ve_p}\varphi_k(\ve_p)$ is weakly convergent in $H^1(O,g)$ and, by Rellich's theorem, norm convergent in $L^2(O,g)$. Hence, $\lim_{p\to\infty}P_{\ve_p}\varphi_k(\ve_p)$ is an eigenfunction for $\Delta_g$ on $H^1(O,g)$ corresponding to eigenvalue $\lim_{p\to\infty}\lambda_k^N(U_\ve, g_{\vert U_{\ve}})$, implying that
$\lim_{p\to\infty}\lambda_k^N(U_\ve, g_{\vert U_{\ve}}) \geq \lambda_k(O,g)$.
\end{proof}


\section{Approximating the spectrum on forms}
\label{sec:SpecDefDim3Form}

The goal of this section is to prove Theorem~\ref{thrm:SpecDefForms}, an analogue of Theorem~\ref{thrm:SpecDefDim3} for the case of the Laplacian acting on differential forms, generalizing \cite[Th\'{e}or\`{e}me 11]{JammesPrecrip} to orbifolds. To begin, recall from Section~\ref{subsec:LaplacianBoundary} that for a closed oriented Riemannian orbifold $(O,g)$ of dimension $n$, $\lambda_{p,k}(O,g)$ denotes the eigenvalue spectrum of the Laplacian acting on $p$-forms and $\mu_{p,k}(O,g)$ denotes the spectrum of the Laplacian on coexact $p$-forms. Recall that the spectrum of the Laplacian on $p$-forms for all $p=0,1,\ldots,n$ is determined by the spectrum on coexact $p$-eigenforms for $p=0,1,\dots, \lfloor (n-1)/2\rfloor$; see Remark~\ref{rem:ExactCoexactSpec}.

Recall \cite[Theorems~1 and 2]{SatakeGenManfld} that de Rham cohomology generalizes readily to orbifolds and is canonically isomorphic to the singular cohomology of the underlying space.
Following \cite{JammesPrecrip}, define the quotient vector space
\[
    H^p(U/O) : = \left\{ \omega \,\vert\, \omega\in\Omega^p(\overline{U}),
    \; d\omega_{\vert\overline{U}} = 0\right\} \big/
    \left\{\omega_{\vert\overline{U}} \,\vert\, \omega\in\Omega^p(O),
    \; d\omega = 0 \mbox{ on $O$}\right\},
\]
where $\overline{U}\subset O$ is a compact orbifold domain with interior $U$ and smooth manifold boundary
$\Theta:=\partial\overline{U}=\partial({O\smallsetminus U})$.
We then have the following generalization of \cite[Th\'{e}or\`{e}me 11]{JammesPrecrip}.

\begin{theorem}
\label{thrm:SpecDefForms}
Let $(O, g)$ be a closed oriented Riemannian orbifold of dimension $n \geq 3$ and let $\overline{U}\subset O$
be a compact orbifold domain with interior $U$ and smooth manifold boundary
$\Theta:=\partial\overline{U}=\partial({O\smallsetminus U})$. Then there exists a sequence of metrics $g_m$ such that
\begin{enumerate}
\item   $\lim_{m\to\infty} \operatorname{Vol}(O, g_m) = \operatorname{Vol}(O,g)$.
\item   Letting $s_p$ denote the dimension of the cohomology group $H^p(U/O)$, we have for
        $p \leq \lfloor(n-3)/2\rfloor$ that
        $\lim_{m\to\infty} \mu_{p,k} (O,g_m) = 0$ when $k\leq s_p$,
        and
        $\lim_{m\to\infty} \mu_{p,k+s_p} (O,g_m) = \mu_{p,k}(O,g)$ when
        $k \geq 1$.
\item   In addition, if for $1\leq p \leq \lfloor (n-3)/2\rfloor$ the $\mu_{p,k}(\overline{U},g)$
        satisfy hypothesis \eqref{eq:HypAst} for some $\delta$, $M$, and $N$, then for every $\ve>0$,
        there exists an $m$ such that
        the $N$-spectral defect of the Laplacians on $O$ and $\overline{U}$ with respect to the metric $g_m$,
        and with Neumann boundary conditions on $\overline{U}$, is $\leq \epsilon$.
\end{enumerate}
\end{theorem}

As stated in \cite[p.973]{JammesPrecrip}, Theorem~\ref{thrm:SpecDefForms} does not hold for $p=\lfloor(n-1)/2\rfloor$.
Note that as in Section~\ref{subsec:SpecDefVerdier}, hypothesis \eqref{eq:HypAst} is satisfied for arbitrarily large $N$.

The Laplacian $\Delta_g$ commutes with both $d$ and $\delta$ and hence sends exact forms to exact forms. Therefore, applying Theorem~\ref{thrm:minmaxNeumannDirichlet} to the case of a closed orbifold so that the boundary conditions are vacuous yields the following; see \cite[Proposition 8]{JammesPrecrip}.

\begin{proposition}
\label{prop:RestrictExactQuadratic}
Let $(O, g)$ be a closed oriented Riemannian orbifold of dimension $n$.
The eigenvalues and eigenspaces of the restriction of the Laplacian to the space of exact forms
are those of the quadratic form
\[
    \tilde{q}(\omega)   =   \Vert\omega\Vert_{L^2\Omega^p(O,g)}^2
\]
relative to the norm
\[
    \Vert\omega\Vert_g    =   \inf\limits_{d\varphi = \omega} \Vert\varphi\Vert_{L^2\Omega^p(O,g)}.
\]
\end{proposition}

We can now proceed with the proof of Theorem~\ref{thrm:SpecDefForms}, which closely follows that of Jammes in \cite[Th\'{e}or\`{e}me 11]{JammesPrecrip} and is similar to the proof of Theorem~\ref{thrm:SpecDefDim3}. We therefore summarize the proof, making note of the places where one must check that the argument continues to hold for orbifolds.

\begin{proof}[Proof of Theorem~\ref{thrm:SpecDefForms}]
As in the proof of \cite[Th\'{e}or\`{e}me 11]{JammesPrecrip}, fix $\ve >0$ and define the metrics $g_\ve$ and $g_{m,\ve}$ on $O$ as in Equation~\eqref{eq:Defgve}.
The three main steps of the proof are the same as those carried out for functions in the proof of Theorem~\ref{thrm:SpecDefDim3} above. The first step is to show that for fixed $\ve$, $g_{\ve}$ can be approximated by a smooth metric $g_{m,\ve}$ so that as $m\to \infty$, the $N$-spectral defect between the quadratic forms associated to $g_{\ve}$ and $g_{m,\ve}$ can be made arbitrarily small.  The proof of this step for orbifolds is identical to that given in \cite{JammesPrecrip}.

The second and third steps of the proof require that we decompose the space $\mathcal{E}^{p+1}(O)$ of exact $(p+1)$-forms into subspaces $\mathcal{K}_0=\mathcal{K}_1\oplus\mathcal{K}_2$ and $\mathcal{K}_{\infty}$.  As explained in \cite{JammesPrecrip}, the decomposition ultimately rests on \cite[Proposition 6]{JammesPrecrip}, which provides the existence of a coexact primitive with minimal $L^2$-norm for every exact form in $\mathcal{E}^{p+1}(O)$.
Because we have shown that the Hodge decomposition theorem holds for orbifolds, Theorem~\ref{thrm:Hodge}, the proof of \cite[Proposition 6]{JammesPrecrip} carries directly to the orbifold setting,
and the classes in $H^p(U/O)$ can be represented by forms that are harmonic on $U$.
Thus, the decomposition of $\mathcal{E}^{p+1}(O)$ into $\mathcal{K}_{0}\oplus\mathcal{K}_{\infty}$ with $\mathcal{K}_0=\mathcal{K}_1\oplus\mathcal{K}_2$ can be carried out for orbifolds exactly as for the manifold case as explained in Jammes.

From here, the second step uses Theorem~\ref{thrm:1.7-Verdier-first} to show that for small $\ve$, $Q$ and $Q_{\vert \mathcal{K}_0}$ have $N$-spectral defect arbitrarily small.  The third step then uses Theorem~\ref{thrm:Thm-1.8-Verdier} to show that $Q_{\vert \mathcal{K}_0}$ and $g_{\vert \overline{U}}$ have $N$-spectral defect arbitrarily small.  With the decomposition of $\mathcal{E}^{p+1}(O)$ into $\mathcal{K}_{0}\oplus\mathcal{K}_{\infty}$ established, the proofs of the second and third steps carry verbatim from \cite{JammesPrecrip}.  Finally, by an application of Lemma~\ref{lem:NSpecDefTransitive}, we conclude that the statement of Theorem~\ref{thrm:SpecDefForms} holds.
\end{proof}


\section{Collapsing connected sums}
\label{sec:ConnectedSums}

We now apply Theorems~\ref{thrm:SpecDefDim3} and \ref{thrm:RauchTaylor} in the case of functions and Theorem~\ref{thrm:SpecDefForms} in the case of forms to create sequences of orbifolds whose Laplace spectra approach that of a manifold, and vice versa. We begin with the following construction.

Let $(O_i,g_i)$, $i=1,2$, be closed oriented Riemannian orbifolds, both of dimension $n\geq 3$.  After possibly scaling $g_1$ or $g_2$, we may assume for both $(O_1,g_1)$ and $(O_2,g_2)$ that the injectivity radius is greater than $2$. For a point $x_i\in O_i$, denote by $B(x_i,r)$ the open ball of radius $r$ about $x_i\in O_i$. Choose a point $p_1\in O_1$ such that $p_1$ is either nonsingular or an isolated singular point. Choose a point $p_2$ in $O_2$ such that the boundary of $B(p_2,1)$ contains no singular points. We assume for simplicity that within a neighborhood of $B(p_i,2)$, the metric $g_i$ is Euclidean (see \cite[p.~548]{annecolbois}). Again by rescaling, we may further assume that for $i=1,2$, the ball $B(p_i,2)$ is completely contained within a single orbifold chart centered at $p_i$ for $i=1,2$.

We now construct the connected sum $O_\ve$ of $O_1$ and $O_2$. Let $O_i(r)=O_i\smallsetminus B(p_i,r)$ and suppose that $\ve<1$.  If $S^{n-1}(r)\subset \mathbb{R}^n$ denotes the sphere of radius $r$ with the standard metric $h_r$ inherited from the Euclidean metric on $\mathbb{R}^n$, then letting $\partial O_i(r)$ denote the boundary of $O_i(r)$ with inherited boundary metric $\partial g_i$, by our hypothesis that $g_i$ is Euclidean on $B(p_i,2)$, we see that $(\partial O_1(\ve), \partial g_1)$ is isometric to $(S^{n-1}(\ve), h_{\ve})$ and $(\partial O_2(1), \ve^2 \partial g_2)$ is isometric to $(S^{n-1}(1), \ve^2 h_1)$.  Furthermore, $(S^{n-1}(\ve), h_{\ve})$ can be mapped isometrically to $(S^{n-1}(1), \ve^2 h_1)$ via the restriction of the map $\mathbb{R}^n\to\mathbb{R}^n$ given by $x\mapsto \ve^{-1}x$.  Taking the composition of these maps, $\varphi_{\ve}:(\partial O_1(\ve), \partial g_1)\to (\partial O_2(1), \ve^2 \partial g_2)$, allows us to create the connected sum of $(O_1, g_1)$ and $(O_2, \ve^2g_2)$ by excising balls about $p_1$ and $p_2$ and identifying the boundaries of $(O_1(\ve), g_1)$ and $(O_2(1),\ve^2g_2)$ via $\varphi_{\ve}$ so that the resulting sum is oriented.  Call this connected sum $(O_{\ve}, g_{\ve})$, where
\[
g_{\ve}=
\begin{cases}
g_1 & \text{ on $O_1(\ve)$,}\\
\ve^2 g_2 &\text{ on $O_2(1)$}.
\end{cases}
\]

Under this construction, $(O_\ve, g_\ve)$ is a closed, smooth, oriented orbifold, and the metric $g_{\ve}$ is a connected
sum metric in the sense of Definition~\ref{def:ConnSum} that is not differentiable along the glued boundary. Hence, the spectrum
of the Laplacian $\Delta_{\ve}$ associated to $g_{\ve}$ satisfies the hypotheses of Theorem~\ref{thrm:spectraltheorem}.

After using this construction in Sections~\ref{subsec:CollapsFunctions1} through \ref{subsec:CollapsForms} to produce sequences of orbifolds with singular points (resp. manifolds) whose Laplace spectra converge to that of a manifold (resp. orbifold with singular points), in Section~\ref{subsec:CollapsPrescrib} we use Theorem~\ref{thrm:SpecDefDim3} to generalize to orbifolds a result of Colin de Verdi\`ere about prescribing the first $N$ elements of the spectrum.


\subsection{Collapsing for the Laplacian acting on functions}
\label{subsec:CollapsFunctions1}

We now show that with the construction of $(O_{\ve},g_{\ve})$ as above, as $\ve$ tends to zero, there is a family of metrics on $O_{\ve}$ agreeing with $g_1$ on $O_1$ so that the eigenvalues of the Laplacian acting on functions on $O_{\ve}$ tend to those of $O_1$ with its original metric.

\begin{theorem}
\label{thrm:RTsmoothconvergence}
Let $(O_i,g_i)$, $i=1,2$, be closed oriented Riemannian orbifolds of dimension $n\geq 3$ and for each $\ve > 0$, let $(O_{\ve}, g_{\ve})$ denote the connected sum as described above. Let $\eta > 0$ and let $N>0$ be an integer. There exist $\ve>0$ and a smooth metric $h_{\eta,\ve,N}$ on $O_{\ve}$ with $h_{\eta,\ve,N \vert O_1(\ve)} = g_{1\vert O_1(\ve)}$ and such that for all $k = 0,1,\dots, N$,
\[
    \vert \lambda_k (O_{\ve}, h_{\eta,\ve,N}) - \lambda_k (O_1,g_1)\vert < \eta.
\]
\end{theorem}
\begin{proof}
By Theorem~\ref{thrm:RauchTaylor}, we have for each $k \geq 0$ that
$\lim_{\ve\to 0}\lambda_k^N (O_1(\ve), g_{1\vert O_1(\ve)}) = \lambda_k(O_1, g_1)$ where we recall that
$\lambda_k^N (O_1(\ve), g_{1\vert O_1(\ve)})$ denotes the $k$th eigenvalue of the Neumann problem on
$(O_1(\ve), g_{1\vert O_1(\ve)})$ and $\lambda_k(O_1, g_1)$ denotes the $k$th eigenvalue of the Laplacian on $(O_1, g_1)$.
Hence, for fixed $\eta > 0$ and $N > 0$,
there exists $\ve>0$ such that the $N$-spectral defect of the Neumann problem on $(O_1(\ve),g_{1\vert O_1(\ve)})$ and the Laplacian on $(O_1,g_1)$ is less than $\eta/2$.

As $\lim_{k\to\infty}\lambda_k(O_1, g_1)=\infty$ by \cite[Proposition 3.2]{chiang}, there are $M$, $N^\prime \geq N$, and $\delta > 0$
such that hypothesis \eqref{eq:HypAst} holds for the $\lambda_k(O_1, g_1)$. Then by decreasing $\ve$ if necessary,
hypothesis \eqref{eq:HypAst} holds as well for the $\lambda_k^N (O_1(\ve), g_{1\vert O_1(\ve)})$ (with $2M$, $N^\prime$,
and $\delta/2$). Then there exists by Theorem~\ref{thrm:SpecDefDim3} a metric $h_{\eta,\ve,N}$ on $O_\ve$
with $h_{\eta,\ve,N\vert O_1(\ve)} = g_{1\vert O_1(\ve)}$ and such that the $N$-spectral defect of the Laplacian on $(O_\ve, h_{\eta,\ve,N})$ and of the Neumann problem on $(O_1(\ve),g_{1\vert O_1(\ve)})$ is less than $\eta/2$.  Thus, it follows that the $N$-spectral defect of the Laplacian on $(O_\ve, h_{\eta,\ve,N})$ and the Laplacian on $(O_1,g_1)$ is less than $\eta$ as desired.
\end{proof}

As a special case of Theorem~\ref{thrm:RTsmoothconvergence}, we prove the existence of a sequence of orbifolds with singular points whose Laplace spectra converge to that of a smooth manifold.

\begin{theorem}
\label{thrm:RTorbifoldstomanifold}
For each $N > 0$ and $m \geq 2$, there is a sequence $\{(O^\ell, h^\ell)\}_{\ell\in\mathbb{N}}$ of $2m$-dimensional closed oriented orbifolds with singular points and a closed oriented manifold $(M,g)$ such that as $\ell\to \infty$, the Laplace spectra of $(O^\ell,h^\ell)$ converge to that of $(M, g)$ in the sense that for all $k=0,1,\dots, N$,
\[
    \vert \lambda_k(O^\ell,h^\ell) - \lambda_k(M,g)\vert< 1/\ell.
\]
\end{theorem}
\begin{proof}
In the construction above, let $(O_1, g_1) = (M,g)$ where $(M,g)$ is a closed oriented Riemannian manifold and take $(O_2, g_2)$ to be an orbifold with singular points such that the singular set of $O_2$ lies outside $B(p_2,1)$.  Let $\eta=1/\ell$.  By Theorem~\ref{thrm:RTsmoothconvergence}, there exist $\ve_{1/\ell}>0$ and metric $h_{1/\ell, \ve_{1/\ell}, N}$ on $O_{\ve_{1/\ell}}$ such that for all $k = 0,1,\dots, N$,
\[
    \vert \lambda _k(O_{\ve_{1/\ell}}, h_{1/\ell,\ve_{1/\ell},N}) - \lambda_k(M,g)\vert < 1/\ell.
\]
By shrinking if necessary, we may assume $\ve_{1/(\ell+1)} < \ve_{1/\ell} < 1/\ell$ for each $\ell$.
Then take $(O^{\ell},h^{\ell}) = (O_{\ve_{1/\ell}}, h_{1/\ell, \ve_{1/\ell},N})$.
\end{proof}

Finally, let us note that for every $m \geq 1$, there is a compact $2m$-dimensional orbifold with a single singular
point. When $m = 1$, the well-known teardrop is such an example. To see this for $m \geq 2$, let $s, a_1, \ldots, a_t$
be positive integers such that $s \geq 2$ and $\gcd(s,a_i)=1$ for each $i$. Recall that the \emph{lens space}
$L(t,a_1\ldots,a_t)$ is the quotient of $S^{2t+1} \subset \C^{t+1}$ by the action of $\Z_s$ generated by
$(z_1,\ldots,z_{t+1})\mapsto(\zeta^{a_1}z_1,\ldots,\zeta^{a_t}z_t, \zeta z_{t+1})$ where $\zeta$ is a primitive
$s$th root of unity. By \cite[Theorems 4.4, 4.10, and 4.16]{SarkarSuhLens}, there exists for each $t\geq 1$ a lens
space $L(t,a_1\ldots,a_t)$ that is the boundary of an oriented manifold $M$ that is compact by construction.
Letting $B(t,a_1\ldots,a_t)$ denote the orbifold given by the quotient of the closed unit ball
$\mathbb{B}^{2t+2} \subset \C^{t+1}$ by the same action and noting that $B(t,a_1\ldots,a_t)$ has a single
singular point with isotropy $\Z_s$, we have $\partial B(t,a_1\ldots,a_t) = L(t,a_1\ldots,a_t)$. Hence,
identifying $\partial M$ with $\partial B(t,a_1\ldots,a_t)$ yields a closed, oriented orbifold with a single
singular point as claimed.

With this, suppose that $(O_1,g_1)$ is an oriented orbifold with a single singular point $p_1$.
If we excise a neighborhood about $p_1$ and take the connected sum with a manifold $(M,h)$ in place of
$(O_2,g_2)$, an argument identical to the proof of Theorem~\ref{thrm:RTorbifoldstomanifold}
yields the following.

\begin{theorem}
\label{thrm:RTmanifoldstoorbifold}
For each $N > 0$ and $m \geq 2$, there is a sequence $\{(M^\ell,h^\ell)\}_{\ell\in\mathbb{N}}$ of $2m$-dimensional closed oriented manifolds and a closed oriented orbifold $(O,g)$ with a single singular point such that as $\ell\to\infty$, the spectra of $(M^\ell,h^\ell)$ converge to that of $(O, g)$ in the sense that for all $k=0,1,\dots, N$,
\[
    \vert \lambda_k(M^\ell,h^\ell) - \lambda_k(O,g)\vert < 1/\ell.
\]
\end{theorem}


\subsection{Another approach to collapsing for functions}
\label{subsec:CollapsFunctions2}

As in Section~\ref{subsec:CollapsFunctions1}, we consider again the relationship between the eigenvalues of the Laplacian acting on functions on $O_{\ve}$ and on $O_1$.  If we relax the condition that the metric on $O_1$ be the original metric $g_1$, we achieve the following using only Theorem~\ref{thrm:SpecDefDim3}.

\begin{theorem}
\label{thrm:smoothconvergence}
Let $(O_i, g_i)$, $i=1,2$, be closed oriented Riemannian orbifolds of dimension $n\geq 3$ and let $(O_{\ve}, g_{\ve})$ denote the connected sum as described above.
Let $\ve>0$, $\eta>0$, and let $N > 0$ be an integer. There is a smooth metric $h_{\eta,\ve,N}$ on the connected sum $O_{\ve}$ and a smooth metric $g_{\eta,\ve,N}$ on $O_1$ such that $h_{\eta,\ve,N \vert O_1(\ve)} = g_{\eta,\ve,N \vert O_1(\ve)} = g_{1\vert O_1(\ve)}$, i.e., all three metrics coincide on $O_1(\ve)$, and such that for all $k=0, 1,\dots, N$,
\[
    \vert \lambda_k(O_{\ve},h_{\eta,\ve,N})-\lambda_k(O_1,g_{\eta,\ve,N})\vert <\eta.
\]
\end{theorem}
\begin{proof}
Considering $O_1(\ve)$ as a subset of the smooth orbifold $O_{\ve}$, extend the metric $g_{1\vert O_1(\ve)}$ arbitrarily
to a smooth metric $\tilde{g_{\ve}}$ on $O_{\ve}$. As in the proof of Theorem~\ref{thrm:RTsmoothconvergence},
as $\lim_{k\to\infty}\lambda_k(O_{\ve}, \tilde{g_{\ve}})=\infty$, by increasing $N$ if necessary, there is a $\delta$ and $M$
such that hypothesis \eqref{eq:HypAst} holds for the $\lambda_k(O_{\ve}, \tilde{g_{\ve}})$, hence for the
$\lambda_k^N(O_1(\ve),g_{1\vert O_1(\ve)})$. Applying Theorem~\ref{thrm:SpecDefDim3} to the smooth metric $\tilde{g_{\ve}}$ on $O_{\ve}$ with $U = O_1(\ve) \subset O_{\ve}$, there is a smooth metric $h_{\eta,\ve,N}$ on $O_{\ve}$ that agrees with $g_1$ on $O_1(\ve)$ such that the $N$-spectral defect of the Neumann problem on $(O_1(\ve),g_{1\vert O_1(\ve)})$ and of the Laplacian on $(O_{\ve}, h_{\eta,\ve,N})$ is less than $\eta/2$. Similarly, applying Theorem~\ref{thrm:SpecDefDim3} to the smooth metric $g_1$ on $O_1$ with $U = O_1(\ve) \subset O_1$, there is a smooth metric $g_{\eta,\ve,N}$ on $O_1$ that agrees with $g_1$ on $O_1(\epsilon)$ such that the $N$-spectral defect of the Neumann problem on $(O_1(\ve),g_{1\vert O_1(\ve)})$ and $(O_1,g_{\eta,\ve,N})$ is less than $\eta/2$. It follows that the $N$-spectral defect of the Laplacian on $(O_{\ve}, h_{\eta,\ve,N})$ and the Laplacian on $(O_1, g_{\eta,\ve,N})$ is less than $\eta$, completing the proof.
\end{proof}

Using the same logic as for Theorems~\ref{thrm:RTorbifoldstomanifold} and \ref{thrm:RTmanifoldstoorbifold} together with Theorem~\ref{thrm:smoothconvergence}, we have the following two results.

\begin{theorem}
\label{thrm:orbifoldstomanifold}
For each $N > 0$ and $m \geq 2$, there is a sequence $\{(O^\ell, h^\ell)\}_{\ell\in\mathbb{N}}$ of $2m$-dimensional closed oriented orbifolds with singular points, a closed oriented manifold $M$, and a sequence $g^\ell$ of metrics on $M$ such that as $\ell\to \infty$, the spectra of $(O^\ell,h^\ell)$ and $(M, g^\ell)$ converge in the sense that for all $k=0,1,\dots, N$,
\[
    \vert \lambda_k(O^\ell,h^\ell)-\lambda_k(M,g^\ell)\vert < 1/\ell.
\]
Moreover, the metrics $g^\ell$ converge pointwise to a smooth metric $g$ on the complement of a point in $M$.
\end{theorem}

\begin{theorem}
\label{thrm:manifoldstoorbifold}
For each $N > 0$ and $m \geq 2$, there is a sequence $\{(M^\ell,h^\ell)\}_{\ell\in\mathbb{N}}$ of $2m$-dimensional closed oriented manifolds, a $2m$-dimensional closed oriented orbifold $O$ with a single singular point, and a sequence of metrics $g^\ell$ on $O$ such that as $\ell\to\infty$, the spectra of $(M^\ell,h^\ell)$ and $(O, g^\ell)$ converge in the sense that for all $k=0,1,\cdots, N$,
\[
    \vert \lambda_k(M^\ell,h^\ell)-\lambda_k(O,g^\ell)\vert < 1/\ell.
\]
Moreover, the metrics $g^\ell$ converge pointwise to a metric $g$ on the complement of the single singular point in $O$.
\end{theorem}


\subsection{Collapsing for the Laplacian acting on forms}
\label{subsec:CollapsForms}

In this section, we observe that using Theorem~\ref{thrm:SpecDefForms}, we may extend
the results of Section~\ref{subsec:CollapsFunctions2} to the spectrum of the Laplacian
acting on forms. Specifically, we have the following.

\begin{theorem}
\label{thrm:smoothconvergencePForms}
Let $(O_i, g_i)$, $i=1,2$, be closed oriented Riemannian orbifolds of dimension $n\geq 3$ and for each $\ve>0$, let $(O_{\ve}, g_{\ve})$ denote the connected sum as described above. For $1 \leq p \leq \lfloor(n-3)/2\rfloor$ and each $\ve>0$, $\eta>0$, and integer $N > 0$, there is a smooth metric $h_{\eta,\ve,N}$ on the connected sum $O_{\ve}$ and a smooth metric $g_{\eta,\ve,N}$ on $O_1$ such that $h_{\eta,\ve,N \vert O_1(\ve)} = g_{\eta,\ve,N \vert O_1(\ve)} = g_{1\vert O_1(\ve)}$, i.e., all three metrics coincide on $O_1(\ve)$, and such that for all $k = 0, 1,\dots, N$,
\[
    \vert \lambda_{p,k}(O_{\ve},h_{\eta,\ve,N}) - \lambda_{p,k}(O_1,g_{\eta,\ve,N})\vert < \eta.
\]
\end{theorem}

The proof is identical to that of Theorem~\ref{thrm:smoothconvergence}, except that we apply Theorem~\ref{thrm:SpecDefForms} in place of Theorem~\ref{thrm:SpecDefDim3}. Note that by Theorem~\ref{thrm:spectraltheoremNeumannDirichlet}, $\lim_{k\to\infty}\lambda_{p,k}^N(O_1(\ve),g_{1\vert O_1(\ve)})=\infty$ so that by increasing $N$ if necessary, there is a $\delta>0$ and $M$ such that the $\lambda_{p,k}^N(O_1(\ve),g_{1\vert O_1(\ve)})$ satisfy hypothesis \eqref{eq:HypAst}. With this, we can again apply the methods used to prove Theorems~\ref{thrm:orbifoldstomanifold} and \ref{thrm:manifoldstoorbifold} to obtain the following.

\begin{theorem}
\label{thrm:orbifoldstomanifoldForms}
For each $N > 0$ and $m \geq 2$, there is a sequence $\{(O^\ell, h^\ell)\}_{\ell\in\mathbb{N}}$ of $2m$-dimensional closed oriented orbifolds with singular points, a closed oriented manifold $M$, and a sequence $g^\ell$ of metrics on $M$ such that as $\ell\to \infty$, the spectra of the Laplacian acting on $p$-forms, $1 \leq p \leq \lfloor(n-3)/2\rfloor$, on
$(O^\ell,h^\ell)$ and $(M, g^\ell)$ converge in the sense that for all $k=0,1,\dots, N$,
\[
\vert \lambda_{p,k}(O^\ell,h^\ell)-\lambda_{p,k}(M,g^\ell)\vert< 1/\ell.
\]
Moreover, the metrics $g^\ell$ converge pointwise to a smooth metric $g$ on the complement of a point in $M$.
\end{theorem}

\begin{theorem}
\label{thrm:manifoldstoorbifoldForms}
For each $N > 0$ and $m \geq 2$, there is a sequence $\{(M^\ell,h^\ell)\}_{\ell\in\mathbb{N}}$ of $2m$-dimensional closed oriented manifolds, a $2m$-dimensional closed oriented orbifold $O$ with a single singular point, and a sequence of metrics $g^\ell$ on $O$ such that as $\ell\to\infty$, the spectra of the Laplacian acting on $p$-forms, $1 \leq p \leq \lfloor(n-3)/2\rfloor$,
on $(M^\ell,h^\ell)$ and $(O, g^\ell)$ converge in the sense that for all $k=0,1,\cdots, N$,
\[
    \vert \lambda_{p,k}(M^\ell,h^\ell) - \lambda_{p,k}(O,g^\ell)\vert < 1/\ell.
\]
Moreover, the metrics $g^\ell$ converge pointwise to a metric $g$ on the complement of the singular point in $O$.
\end{theorem}


\subsection{Prescribing a finite part of the spectrum via collapsing}
\label{subsec:CollapsPrescrib}

In this section, we briefly discuss the question of prescribing the first $N$ nonzero eigenvalues of the Laplacian of a closed oriented orbifold. We appreciate the suggestion from an anonymous referee that these results be included here.

The main result of \cite{VerdierMultiplicite} is that for any closed manifold $M$ of dimension $n\geq 3$ and any positive integer $N$, there is a Riemannian metric on $M$ such that the first nonzero eigenvalue of the Laplacian of $(M,g)$ acting on functions has multiplicity $N$. This was later generalized to show that the first $N$ eigenvalues of $M$ can be prescribed to be any positive values, see \cite{VerdierConstrucDonee}, and was generalized in \cite{JammesPrecrip} to the case of the Laplacian acting on $p$-forms for $1 \leq p < n/2$ when $n \geq 6$.

Here, we indicate how these results can be extended to the case of closed orbifolds. For simplicity, we restrict to the case of the Laplacian acting on functions. Succinctly, Colin de Verdi\`ere's approach involves embedding a graph or surface with the desired spectrum into the manifold, extending the metric globally while maintaining a small $N$-spectral defect, and then using the weak Arnold's hypothesis (see below) to show that the metric on the graph or surface can be perturbed so that the spectrum of the resulting Laplacian on $M$ obtains the desired values. This approach extends directly to the orbifold case by using Theorem~\ref{thrm:SpecDefDim3} and performing the embedding away from the singular set.

Choose $0 < \lambda_1 \leq \lambda_2 \leq\cdots\leq \lambda_N$.
By \cite[Section 2]{VerdierConstrucDonee}, there is a compact neighborhood $B$ of $b_0 := (\lambda_1,\ldots,\lambda_N) \in \R^N$,
a closed surface $S$, and a family of metrics $(g_b)_{b\in B}$ on $S$ satisfying the \emph{weak Arnold's hypothesis (WAH)}.
Specifically, let $E_b$ denote the $b$-eigenspace of the Laplacian $\Delta_{b}$ associated to the metric $g_{b}$, let $\mathcal{Q}(E_b)$ denote the space of real quadratic forms on $E_b$, and let $\Phi\co B\to\mathcal{Q}(E_{b_0})$ be the map that assigns to each $b\in B$ the quadratic form on $E_{b_0}$ associated to $\Delta_b$ via the natural isometry of $E_b$ onto $E_{b_0}$. Then $\Phi(b_0)$ satisfies WAH relative to the family $(\Delta_b)_{b\in B}$ if $\Phi$ is essential at $\Phi(b_0)$, meaning that there is an $\eta > 0$ such that for all continuous $\Psi\co B\to\mathcal{Q}(E_{b_0})$ with $\|\Psi-\Phi\|_{L^\infty(B)} < \eta$, $\Phi(b_0)$ is in the image of $\Psi$. See \cite[Section 1]{VerdierArnold} for more details.

Now, let $O$ be a closed oriented orbifold of dimension $n\geq 3$. By embedding $S$ into an orbifold chart at a nonsingular point of $O$, we may define an embedding $\iota\co S \to O$ such that a tubular neighborhood $U$ of $\iota(S)$ is an orbifold domain that does not intersect the singular set of $O$. For each $b\in B$, define a metric $\tilde{g}_b$ on $U$ as a product metric with $g_b$ on $\iota(S)$ so that the first $N$ nonzero eigenvalues of the Neumann problem on $\overline{U}$ are given by $b$. Then by Theorem~\ref{thrm:SpecDefDim3}, there exists for each $b\in B$ a smooth metric $g_b^O$ on $O$ that extends $\tilde{g}_b$ and such that the $N$-spectral defect of the Neumann problem on $\overline{U}$ and the Laplacian on $(O, g_b^O)$ is $\leq\eta/2$. Moreover, following the proof of Theorem~\ref{thrm:SpecDefDim3}, we may by shrinking $B$ if necessary choose the $g_\ve$ and approximating functions $F_m^\ve$ uniformly for each $b\in B$ so that the $g_b^O$ coincide on $O\smallsetminus U$. Let $\Psi\co B\to \mathcal{Q}(E_{b_0})$ assign to each $b\in B$ the quadratic form associated to the Laplacian of $(O, g_b^O)$ acting on the eigenspace of the first $N$ nonzero eigenvalues, which we identify with $E_{b_0}$ via the natural isometry, and then $\Psi$ is continuous. Hence $\|\Psi-\Phi\|_{L^\infty(B)} < \eta$, which implies that there is a $b\in B$ such that $\Psi(b) = \Phi(b_0)$. That is, the first $N$ nonzero eigenvalues of the Laplacian on $(O, g_b^O)$ are given by $b_0 = (\lambda_1,\ldots,\lambda_N)$. We summarize this observation with the following.

\begin{theorem}
\label{thrm:PrescribNSpec}
Let $0 < \lambda_1 \leq \lambda_2 \leq\cdots\leq \lambda_N$ be real numbers, and let $O$ be a closed oriented orbifold of dimension $n\geq 3$. Then there exists a Riemannian metric $g$ on $O$ such that the first $N$ nonzero eigenvalues of the Laplacian on $(O,g)$ are given by the $\lambda_1,\ldots,\lambda_N$.
\end{theorem}


\bibliographystyle{amsplain}

\begin{thebibliography}{14}

\bibitem{alr}
A. Adem, J. Leida, and Y. Ruan, \emph{Orbifolds and stringy topology}.  Cambridge Tracts in Mathematics, 171.  Cambridge University Press, Cambridge, 2007.


\bibitem{anneperturbNeumann}
C. Ann\'e, \emph{Perturbation du spectre {$X\smallsetminus TUB^\epsilon Y$} (conditions de {N}eumann)}.
S\'{e}minaire de {T}h\'{e}orie {S}pectrale et {G}\'{e}om\'{e}trie, {N}o. 4, {A}nn\'{e}e 1985--1986,
Univ. Grenoble I, Saint-Martin-d'H\`eres (1986), 17--23.


\bibitem{annespectre}
C. Ann\'e, \emph{Spectre du laplacien et \'ecrasement d'anses}.  Ann. Sci. \'Ecole Norm. Sup. (4) \textbf{20} (1987), 271--280.


\bibitem{annecolboisoperateur}
C. Ann\'e and B. Colbois, \emph{Op\'erateur de Hodge-Laplace sur des vari\'et\'es compactes priv\'ees d'un nombre fini de boules}.  J. Funct. Anal. \textbf{115} (1993),  190--211.

\bibitem{annecolbois}
C. Ann\'{e} and B. Colbois, \emph{Spectre du Laplacien agissant sur let $p$-formes diff\'erentielles et \'ecrasement d'anses}.  Math. Ann. \textbf{303} (1995),  545--573.

\bibitem{Anne-Post}
C. Ann\'{e} and O. Post, \emph{Wildly perturbed manifolds: norm resolvent
	and spectral convergence}, preprint, \texttt{arXiv:1802.01124 [math.SP]}, 2018.

\bibitem{annetakahashi}
C. Ann\'{e} and J. Takahashi, \emph{$p$-spectrum and collapsing of connected sums}. Trans. Amer. Math. Soc. \textbf{364} (2012),  1711--1735.


\bibitem{AMDGHRS}
T. Arias-Marco, E. B. Dryden, C. S. Gordon, A. Hassannezhad, A. Ray, and E. Stanhope,
\emph{Spectral geometry of the Steklov problem on orbifolds}. Int. Math. Res. Not. IMRN 2019, 90--139.


\bibitem{bailey}
W. Baily,  Jr.,
\emph{The decomposition theorem for V-manifolds}.
Amer. J. Math. \textbf{78} (1956), 862--888.


\bibitem{berlinegetzler}
N. Berline, E. Getzler, and M. Vergne,
\emph{Heat kernels and Dirac operators}. Corrected reprint of the 1992 original. Grundlehren Text Editions. Springer-Verlag, Berlin, 2004.

\bibitem{Bucicovschi}
B. Bucicovschi, \emph{Seeley's theory of pseudodifferential operators on orbifolds},
preprint, 	\texttt{arXiv:math/9912228 [math.DG]}, 2008.


\bibitem{chavel-feldman-domains}
I. Chavel and E A. Feldman, \emph{Spectra of domains in compact manifolds}.   J. Func. Anal. \textbf{30}  (1978), 198--222.

\bibitem{chavel-feldman-less}
I. Chavel and E A. Feldman, \emph{Spectra of manifolds less a small domain}.   Duke Math. J.
\textbf{56}  (1988), 399--414.

\bibitem{ChenRuan}
W. Chen and Y. Ruan,
\emph{A new cohomology theory of orbifold}.
Comm. Math. Phys. \textbf{248} (2004),  1--31.

\bibitem{chiang}
Y.-J. Chiang, \emph{Harmonic maps of $V$-manifolds}.  Ann. Global Anal. Geom. \textbf{8} (1990), 315--344.

\bibitem{colboiselsufi}
B. Colbois and A. El Soufi, \emph{Spectrum of the Laplacian with weights}.
Ann. Global Anal. Geom. \textbf{55} (2019),  149--180.


\bibitem{VerdierMultiplicite}
Y. Colin de Verdi\`ere,
\emph{Sur la multiplicit\'e de la premi\`ere valeur propre non nulle du laplacien}.
Comment. Math. Helv. \textbf{61} (1986), 254--270.


\bibitem{VerdierConstrucDonee}
Y. Colin de Verdi\`ere,
\emph{Construction de laplaciens dont une partie finie du spectre est donn\'{e}e}.
Ann. Sci. \'{E}cole Norm. Sup. (4) \textbf{20} (1987), 599--615.

\bibitem{VerdierConstrucDoneeAvecMult}
Y. Colin de Verdi\`ere,
\emph{Construction de laplaciens dont une partie finie (avec multiplicités) du spectre est donn\'{e}e}.
S\'{e}minaire sur les \'{e}quations aux d\'{e}riv\'{e}es partielles 1986--1987, Exp. No. VII, 6 pp.,
\'{E}cole Polytech., Palaiseau, 1987.

\bibitem{VerdierArnold}
Y. Colin de Verdi\`ere,
\emph{Sur une hypoth\`{e}se de transversalit\'{e} d'Arnolʹd}.
Comment. Math. Helv. \textbf{63} (1988), 184--193.

\bibitem{davies}
E. B. Davies, \emph{Spectral theory and differential operators}.  Cambridge Studies in Advanced Mathematics, 42.
Cambridge University Press, Cambridge, 1995.

\bibitem{doylerossetti}
P. Doyle and J. P. Rossetti, \emph{Isospectral hyperbolic surfaces have matching geodesics}.  New York J. Math. \textbf{14} (2008), 193--204.

\bibitem{dggw}
E. Dryden, C. Gordon, S. Greenwald, and D. Webb, \emph{Asymptotic expansion of the heat kernel for orbifolds}. Michigan Math. J., \textbf{56} (2008), 205--238.

\bibitem{drydenstrohmaier}
E. Dryden and A. Strohmaier, \emph{Huber's theorem for hyperbolic orbisurfaces}.  Canad. Math. Bull. \textbf{52} (2009), 66--71.

\bibitem{dodziuk}
J. Dodziuk,
\emph{Eigenvalues of the Laplacian on forms}.
Proc. Amer. Math. Soc. \textbf{85} (1982), no. 3, 437--443.

\bibitem{evans}
L. C. Evans, \emph{Partial differential equations}. Graduate Studies in Mathematics, 19.  American Mathematical Society, Providence, RI, 1998.

\bibitem{farsi}
C. Farsi, \emph{Orbifold spectral theory}.  Rocky Mountain J. Math. \textbf{31} (2001), 215--235.

\bibitem{gordonrossetti}
C. Gordon and J. P. Rossetti, \emph{Boundary volume and length spectra of Riemannian manifolds: what the
middle degree Hodge spectrum doesn't reveal}.  Ann. Inst. Fourier (Grenoble) \textbf{53} (2003),  2297--2314.

\bibitem{gornet-mcgowan}
R. Gornet and J. McGowan,
\emph{Small eigenvalues of the Hodge Laplacian for three-manifolds with pinched negative curvature}.
In: Spectral problems in geometry and arithmetic (Iowa City, IA, 1997), 29--38,
Contemp. Math., 237, Amer. Math. Soc., Providence, RI, 1999.

\bibitem{hepworth}
R. Hepworth,
\emph{Morse inequalities for orbifold cohomology}.
Algebr. Geom. Topol. \textbf{9} (2009),  1105--1175.


\bibitem{hormander}
L.H\"{o}rmander, \emph{The Analysis of Linear Partial Differential Operators III}.
Grundlehren der mathematischen Wissenschaften 274, Springer Verlag, Berlin, 1985.

\bibitem{JammesPrecrip}
P. Jammes, \emph{Prescription de la multiplicit\'{e} des valeurs propres du laplacien de Hodge-de Rham}.
Comment. Math. Helv. \textbf{86} (2011), 967--984.

\bibitem{linowitzmeyer}
B. Linowitz and J. S. Meyer,
\emph{On the isospectral orbifold-manifold problem for nonpositively curved locally symmetric spaces}.
Geom. Dedicata \textbf{188} (2017), 165--169.

\bibitem{mcgowan}
J. McGowan, \emph{The {$p$}-spectrum of the {L}aplacian on compact hyperbolic three manifolds}.
Math. Ann. \textbf{297} (1993),  725--745.

\bibitem{Post}
O. Post, \emph{Boundary pairs associated with quadratic forms}. Math. Nachr. \textbf{289} (2016),  1052--1099.

\bibitem{rauchtaylor}
J. Rauch and M. Taylor, \emph{Potential and scattering theory on wildly perturbed domains}.  J. Funct. Anal. \textbf{18} (1975), 27--59.


\bibitem{rsw}
J. P. Rossetti, D. Schueth, and M. Weilandt, \emph{Isospectral orbifolds with different maximal isotropy orders}. Ann. Global Anal. Geom. \textbf{34} (2008), 351--366.

\bibitem{SarkarSuhLens}
S. Sarkar and D. Y. Suh,
\emph{A new construction of lens spaces}.
Topology Appl. \textbf{240} (2018), 1--20.

\bibitem{SatakeGenManfld}
I. Satake,
\emph{On a generalization of the notion of manifold}.
Proc. Nat. Acad. Sci. U.S.A. \textbf{42} (1956), 359--363.

\bibitem{schwarz}
G. Schwarz, \emph{Hodge decomposition - a method for solving boundary value problems}.  Lecture Notes in Mathematics, 1607.  Springer-Verlag, Berlin, 1995.

\bibitem{sutton}
C. Sutton, \emph{Equivariant isospectrality and Sunada's method}. Arch. Math. (Basel) \textbf{95} (2010),  75--85.

\bibitem{takahashi}
J. Takahashi, \emph{Collapsing of connected sums and the eigenvalues of the Laplacian}. J. Geom. Phys. \textbf{40} (2002), 201--208.

\bibitem{takahashiGapPForms}
J. Takahashi, \emph{On the Gap between the First Eigenvalues of the Laplacian on Functions and $p$-Forms}. Ann. Global Anal. Geom. \textbf{23} (2003),  13--27.

\bibitem{taylor}
M. Taylor, \emph{Partial differential equations I.  Basic theory. Second edition}.  Applied Mathematical Sciences, 115.  Springer, New York, 2011.

\bibitem{Wei}
G. Wei, \emph{Manifolds with a lower Ricci curvature bound}. Surveys in differential geometry. Vol. XI, 203--227,
Int. Press, Somerville, MA, 2007.


\end{thebibliography}

\end{document}